\documentclass[review]{elsarticle}
\usepackage{graphicx, color}
\usepackage{lineno,hyperref}
\modulolinenumbers[5]
\newtheorem{ex}{Example}

\newtheorem{lm}{Lemma}
\newtheorem{defi}{Definition}
\newtheorem{thm}{Theorem}
\newtheorem{rem}{Remark}

\journal{}

\begin{document}

\begin{frontmatter}

\title{Tails and probabilities for $p$-outside values}
\tnotetext[mytitlenote]{\href{https://pavlinakj.wordpress.com/}{The author is grateful to the bilateral projects Bulgaria - Austria, 2016-2019, "Feasible statistical modelling for extremes in ecology and finance", BNSF, Contract number 01/8, 23/08/2017 and WTZ Project No. BG 09/2017.}.}

\author[myfootnote]{Pavlina Jordanova}
\address{Faculty of Mathematics and informatics, Konstantin Preslavsky University of Shumen,\\
Universitetska 115, Shumen, Bulgaria}
\fntext[myfootnote]{Corresponding author: pavlina\_kj@abv.bg, p.jordanova@shu.bg}




\begin{abstract}
 The task for a general and useful classification of the tails of probability distributions still has no satisfactory solution. Due to lack of information outside the range of the data the tails of the distribution should be described via many characteristics. Index of regular variation is a good characteristic, but it puts too many distributions with very different tail behavior in one and the same class. One can consider for example Pareto($\alpha$), Fr$\acute{e}$chet($\alpha$) and Hill-horror($\alpha$) with one and the same fixed parameter $\alpha > 0$. The main disadvantage of VaR, expectiles, and hazard functions, when we speak about the tails of the distribution, is that they depend on the center of the distribution and on the scaling factor. Therefore they are very appropriate for predicting "big losses", but after a right characterization of the distributional type of "the payoff". When analyzing the tail of the observed distribution we need some characteristic which does not depend on the moments because in the most important cases of the heavy-tailed distributions theoretical moments do not exist and the corresponding empirical moments fluctuate too much. In this paper, we show that probabilities for different types of outside values can be very appropriate characteristics of the tails of the observed distribution. They do not depend on increasing affine transformations and do not need the existence of the moments. The idea origins from Tukey's box plots, and allows us to obtain one and the same characteristic of the  tail of the observed distribution within the whole distributional type with respect to all increasing affine transformations. These characteristics answer the question:
\begin{center}
{\it At what extent we can observe "unexpected" values?}
\end{center}
\end{abstract}


\end{frontmatter}

\section{Motivation and history of the problem}
\label{intro}
The task for a general and useful classification of probability distributions with respect to the their tails  seems to be still open.
Embrechts et al. (1997) \cite{EMK} have made a very useful figure about the relations between different subclasses of heavy-tailed distributions in the sense of the infinite moment generating function for all positive arguments. However, this classification puts too many distributions with very different tail behavior in one and the same class. According to this classification for example Pareto($\alpha$), Fr$\acute{e}$chet($\alpha$) and Hill-horror($\alpha$) distributions, with one and the same fixed parameter $\alpha > 0$ belong to one and the same class of distributions with regularly varying tails with parameter $\alpha$. However, the chance to observe "unexpected" value in these three cases is very different, especially for the Hill-Horror distribution. See Figure \ref{fig:peRalphaE}. Comparison of the corresponding hazard rate functions $r_X(x)$, cumulative distribution functions $F_X(x)$ (c.d.fs), and probability density functions $f_X(x)$ (p.d.fs) when $x$ is fixed, and close to the ends of the support of the corresponding distributions is a relatively good approach, but all these characteristics depend on the center of the distribution and the scale parameters. In order to delete these dependencies usually, we normalize the considered random variable(r.v.) with the variance. However, in order to do this, we need existence not only of the first but also of the second moment of the observed distribution. In the most important cases of heavy-tailed distributions, these moments do not exist and this approach is not applicable.  Therefore we need some characteristics which describe separately the left and the right tail of the distribution and do not depend on the moments. Index of regular variation, in cases when it is meaningful, is not enough. These lead us to the idea about usage of quantiles, VaR and expectiles, described e.g. in Daouia et al. (2018) \cite{Daouia2018a} or Marinelli et al. (2007) \cite{MarinelliCarlo2007}. They are good characteristics, but also depend on the center of the distribution and the scaling factor. Therefore they are very appropriate for predicting "big looses" within a fixed family of distributional type, but first, you need a right characterization of the type of the "the profit and loss" distribution. A long critical review about the kurtosis can be seen in Balanda and MacGillary (1988) \cite{BalandaMacGillary1988}. Harter (1959) \cite{Harter} uses sample quasi-ranges in estimating population standard deviation. Mosteller (1946) \cite{Monsteller1946} and Sarhan (1954) \cite{SarhanGeneral} propose estimators of the mean and standard deviation which are functions of order statistics. They give us the idea to work with something related to the quantile spread.

Due to lack of information outside the range of the data the tails of the distribution should be described via many characteristics.
Through the paper, we show that probabilities for different orders of outside values can be appropriate characteristics for solving this task. Their properties outperform the properties of the kurtosis, tail index and hazard function when speaking about classification with respect to the tail of the observed distribution. The main their advantages are that they do not depend on the center and the scaling factor of the distribution, and do not need the existence of the moments. They are useful for answering the question:
\begin{center}
{\it At what extend we should observe "unexpected" values?}
\end{center}
The idea origins from Tukey's box plots (1977) \cite{Tukey1977} and Balanda and MacGillary's (1990) \cite{BalandaMacGillary} spread-spread plot.  However instead of the quantiles here we use probabilities. This allows us to obtain one and the same characteristic of the tail of the observed distribution within all distributional type with respect to increasing affine transformations.

In Section 2 we define and investigate the general properties of probabilities for $p$-outside values. In Section 3 their explicit forms are calculated and plotted for the most popular probability distributions. Section 4 investigates asymptotic properties of empirical left and right $p$-fences, and the estimators of probabilities for left and right $p$-outside values. A result about strong consistency of relative frequency estimator completes that part. 

Different estimators of the exponent of regular variation are proposed in Hill (1975) \cite{Hill}, Pickands (1975)\cite{Pickands} and Deckers-Einmahl-de Haan Dekkers (1989) \cite{Dekkers1989}, Einmahl and Guillou (2008) \cite{EinmahlGuillou}, t-Hill Stehlik et al. (2010) \cite{Stehlik2010}, Pancheva and Jordanova (2012) \cite{JordanovaPancheva,JordanovaMilan2012}, Jordanova et al. (2016) \cite{OurExtremes} among others. The mean of order $p$ generalization of the t-Hill and Hill statistics is introduced by Beran et al. (2014) \cite{Beran}, Caeiro et al. (2016) \cite{CaeiroGomes} and Paulauskas and Vaiciulis (2017) \cite{Paulauskas}. Another approach can be seen in Huisman et al. (2001) \cite{Huisman} who recommend to correct small-sample bias of Hill estimators via weighted averages of its values for different thresholds.
In Sections 5 the previous results are applied and a completely new approach for estimating the parameter of the heaviness of the tail of the observed distribution is demonstrated. Four of the examples consider cumulative distribution functions (c.d.fs.) with regularly varying  right tail. These are Pareto, Fr$\acute{e}$chet, Log-logistic and Hill-horror cases. It is easy to realize that in order to estimate their indexes of regular variation working with small samples only distribution sensitive estimators can be useful. The main idea of this section is to show that our approach works also in the case when the c.d.f. of the observed r.v. does not have regularly varying right tail. We depict this result via an example about $H_1$ distribution (\ref{H1CDF}) which tail is not regularly varying.  The paper finishes with some conclusive remarks. The proofs are sent to the Appendix section.
All plots and computations are made via software R (2018)\cite{R}.

In this work we do not use the second order regular variation introduced in de Haan and Stadtmueller (1996) \cite{deHaanStadtmueller} for distributions with regularly varying tails, because it is applicable only for huge samples. A very comprehensive study of analyzing extreme values under the second order regularly varying condition can be found e.g. in de Haan and Ferreira (2006) \cite{deHaanFerreira}. The corresponding  properties of the convolutions and the central limit theorem are obtained by Geluk and de Haan (1997) \cite{GelukdeHaan1997}.

Through the paper we use the following notations: $ {\mathop{=}\limits_{}^{d}}$ is for the equality in distribution, $ {\mathop{\to}\limits_{}^{d}}$ means convergence in distribution, ${\mathop{\to}\limits_{}^{a.s.}}$  denotes almost sure convergence.  $\in$ means that the considered r.v. belongs to the corresponding probability type. $\sim$ is asymptotic equivalence. $B(\alpha, \beta) = \int_0^1 x^{\alpha - 1}(1- x)^{\beta-1} dx$ is for the beta function, and $Beta(\alpha, \beta)$ denotes Beta distribution with parameters $\alpha > 0$ and $\beta > 0$.
 $\mathbf{X} \in Par(\alpha, \delta)$, means that a r.v. $\mathbf{X}$ has Pareto c.d.f.
\begin{equation}\label{Pareto}
 F_{\mathbf{X}}(x) = \left\{
                  \begin{array}{ccc}
                    0 & , & x \leq \delta \\
                    1 - \left(\frac{\delta}{x}\right)^{\alpha} & , & x > \delta
                  \end{array}
                \right..
\end{equation}

More general definitions of Pareto distributions, together with very useful descriptions of the relations between them and the most important other distributions could be seen e.g. in Arnold (2015) \cite{Arnold2015}.

$\mathbf{X} \in Fr(\alpha, c)$, means that a r.v. $\mathbf{X}$ has Fr$\acute{e}$chet c.d.f.
\begin{equation}\label{Frechet}
 F_{\mathbf{X}}(x) = \left\{
                  \begin{array}{ccc}
                    0 & , & x \leq 0 \\
                    e^{-(cx)^{-\alpha}} & , & x > 0
                  \end{array}
                \right..
\end{equation}

$\mathbf{X} \in NegWeibull(\alpha, \sigma,  \mu)$ is an abbreviation of the fact that the r.v.  $\mathbf{X}$ has a Negative Weibull c.d.f.
\begin{equation}\label{Weibull}
F_{\mathbf{X}}(x) =\left\{\begin{array}{ccc}
                   exp\left\{-\left(-\frac{x-\mu}{\sigma}\right)^{\alpha}\right\} & , & x \leq \mu\\
                   1 & ,  & x > \mu
                 \end{array}
\right..
\end{equation}

We consider only absolutely continuous distributions. For $p \in [0, 1]$ the theoretical quantile function of the c.d.f. $F$ is defined as $$F^\leftarrow(p) = inf\{x \in R: F(x) \geq p\} = sup \{x \in R: F(x) \leq p\}.$$

Let $X_1, X_2, ..., X_n$ be a sample of independent observations on a r.v. $X$  with c.d.f. $F$. Here we denote the corresponding order statistics by $X_{(1, n)} \leq X_{(2, n)} \leq ... \leq X_{(n, n)}$. In   Parzen (1979) \cite{Parzen}, Hyndman et al. (1996) \cite{Hyndman}, Langford (2006) \cite{Langford} among others, one can find different definitions of empirical $p$-quantiles, $p \in \left[\frac{1}{n+1}, \frac{n}{n+1}\right]$.  We use the following one $F^\leftarrow_n(p):= X_{([(n+1)p],n)}$, where $[a]$ means the integer part of $a$. {\footnote{As it is noticed in Chu (1957) \cite{Chu}, for large samples, these methods are equivalent because we consider only absolutely continuous distributions.}}


\section{Using probabilities for $p$-outside values for characterising the tails of the observed distribution}
\label{sec:1}
The idea about classification of distributions based on quartiles and box plots comes from Tukey (1977) \cite{Tukey1977}, and recently was reminded by Devore (2015) \cite{Devore} and Jordanova and Petkova (2017) \cite{MoniPoli2017}. Here we generalize this concept and introduce one more parameter in the definition of outside values, which allows the researchers to decide at what extend atypical observations would be called "outside value".

Denote by $R_n(F,p) = R_n(X,p)$
$$ = F^\leftarrow_n(1 - p)  + \frac{1-p}{p}(F^\leftarrow_n(1 - p) - F^\leftarrow_n(p)) = \frac{1}{p}F^\leftarrow_n(1 - p) - \frac{1-p}{p} F^\leftarrow_n(p)$$
and by $L_n(F,p) = L_n(X,p)$
$$ = F^\leftarrow_n(p) - \frac{1-p}{p}(F^\leftarrow_n(1 - p) - F^\leftarrow_n(p)) = \frac{1}{p}F^\leftarrow_n(p) - \frac{1-p}{p} F^\leftarrow_n(1 - p),$$
 correspondingly {\bf empirical p-right-} and {\bf empirical  p-left-fence}.

Their sum is equal to $F^\leftarrow_n(1 - p) + F^\leftarrow_n(p)$. The difference $F^\leftarrow_n(1 - p) - F^\leftarrow_n(p)$ between these quantiles is very well known. It is called {\bf empirical quantile spread(quasi range)}, and is considered e.g. in Gumbel(1944)\cite{Gumbel1944}, Monsteller (1946) \cite{Monsteller1946}, and Balanda and MacGillary (1990) \cite{BalandaMacGillary}.
The meaning of these values comes from the expression
$$p[R_n(F,p) - F^\leftarrow_n(p)] = F^\leftarrow_n(1 - p) - F^\leftarrow_n(p) = p[F^\leftarrow_n(1 - p) - L_n(F,p)].$$


It is clear that analogously to Tukey's box-plot (1977) \cite{Tukey1977}, one can use {\bf empirical box plot of order p}.
Its borders are determined via the values
$$L_n(F,p),\,\, F^\leftarrow_n(p),\,\, F^\leftarrow_n(0.5),\,\, F^\leftarrow_n(1 - p),\,\, R_n(F, p).$$
The most frequently $p = 0,25$. This case is partially investigated in the supplementary material of Soza et al. (2019) \cite{JordanovaStehlik2018}.

Sample right  or left $p$-outside values are the observations which fall outside the interval  $[L_n(F,p), R_n(F,p)]$. Their absolute frequencies, strongly depend on the sample size.
\begin{defi}\label{def1} Assume $p \in (0, 0.5]$. We call on observation $Y$ sample(empirical)
\begin{itemize}
\item {\bf right $p$-outside values} if $Y > R_n(X,p);$
\item {\bf left $p$-outside values} if $Y < L_n(X,p).$
\end{itemize}
\end{defi}

For $p = 0.25$, the definition coincides with the one in Devore (2015) \cite{Devore}. He calls them "extreme right" and "extreme left outliers".

Denote by $n_{R}(p, n)$, and $n_{L}(p, n)$ the numbers of these outside values in a sample of $n$ independent observations. According to the A.Kolmogorov's Zero-one law, never mind how small, but strictly positive is $p$ one can almost sure observe such outside values in a large enough sample of observations on a r.v. with many light tailed distributions,  e.g. Gaussian. Therefore the number of outside values is not too informative.
In order to classify distributions with respect to their tail behaviour we propose to compare their theoretical probabilities an observation to be an outside value of the considered type. Denote by
\begin{eqnarray*}
  p_{L,p}(X) &=& p_{L,p}(F) = P(X < L(X,p)),\\
  p_{R, p}(X) &=& p_{R, p}(F) = P(X > R(X,p))
\end{eqnarray*}
the observed r.v. $X$ to be left or right $p$-outside value.
Here, analogously to $R_n(F,p)$  and $L_n(F,p)$ we have denoted by
\begin{eqnarray*}
  R(F,p) &=& R(X,p) = F^\leftarrow(1 - p)  + \frac{1-p}{p}[F^\leftarrow(1 - p) - F^\leftarrow(p)]\\
   &=& \frac{1}{p}F^\leftarrow(1 - p) - \frac{1-p}{p} F^\leftarrow(p)
\end{eqnarray*}
{\bf theoretical p-right fence}  and by
\begin{eqnarray*}
  L(F,p) &=& L(X,p) = F^\leftarrow(p) - \frac{1-p}{p}[F^\leftarrow(1 - p) - F^\leftarrow(p)]\\
   &=& \frac{1}{p}F^\leftarrow(p) - \frac{1-p}{p} F^\leftarrow(1 - p)
   \end{eqnarray*}
{\bf theoretical p-left fence}.
Their properties are analogous to the properties of empirical p-right- and p-left-fences.

Let us note that for any absolutely continuous c.d.f. $F$
$$R(X,0.5) = L(X,0.5) = F^\leftarrow(0.5),$$
is the median and $p_{R, 0.5}(X) = p_{L, 0.5}(X) = 0.5.$

It is not difficult to check that $R(F, p)$ and $L(F, p)$ are monotone. Therefore by monotonicity of probability measures, the characteristics $p_{L,p}(X)$ and $p_{R,p}(X)$ are also monotone.
\begin{thm}\label{thm:monotonicity} For a fixed $p \in (0, 0.5]$ if there exist $f(F^\leftarrow(p)) \in (0, \infty)$, and $f(F^\leftarrow(1 - p)) \in (0, \infty)$, then,
\begin{description}
\item[a)] $L(F, p)$ is increasing in $p$;
\item[b)] $R(F, p)$ is decreasing in $p$;
\item[c)] $p_{L,p}(X)$ and $p_{R,p}(X)$ are non-decreasing in $p$.
\end{description}
\end{thm}

For $p = 0.25$, $p_{L, p}(X)$ and $p_{R, p}(X)$ are correspondingly the probabilities an observation to be left- or right- extreme outlier.  Jordanova and Petkova (2018) \cite{MoniPoli2018} and Soza et al.(2019) \cite{JordanovaStehlik2018} denote these probabilities by $p_{eL} = p_{L, 0.25}$ and $p_{eR} = p_{R, 0.25}$. In that case, the authors  obtain them for Pareto, Fr$\acute{e}$chet, $H_1$, $H_2$, and Hill-Horror distributions. Further on, we generalize these results.

\begin{thm}\label{thm:thm1} Assume $p \in (0, 0.5]$, and $g(x)$ and $F := F_X$ are strictly monotone, well defined, and continuous function in the considered values. The characteristics $p_{L,p}(X)$, $p_{R,p}(X)$,  possess the following properties:
\begin{description}
\item[a)] $p_{L,p}(X) \in [0, p]$ and  $p_{R,p}(X) \in [0, p]$.
\item[b)] $p_{L,p}(X) = p_{L,p}(X + c)$, $p_{R,p}(X) = p_{R,p}(X + c)$, $c \in R$.
\item[c)] If $c > 0$, then $p_{L,p}(cX) = p_{L,p}(X)$, $p_{R,p}(cX) = p_{R,p}(X)$.
\item[d)] If $c < 0$, then $p_{L,p}(cX) = p_{R,p}(X)$, $p_{R,p}(cX) = p_{L,p}(X)$.
\item[e)] If $g(x)$ is continuous and strictly increasing
\begin{eqnarray}
\label{eL}  p_{L,p}(g(X)) &=& F \left\{g^\leftarrow\left[\frac{1}{p} g[F^\leftarrow(p)] - \frac{1-p}{p}g[F^\leftarrow(1 - p)]\right]\right\}, \\
\label{eR}  p_{R,p}(g(X)) &=& 1 - F \left\{g^\leftarrow\left[\frac{1}{p} g[F^\leftarrow(1-p)] - \frac{1-p}{p}g[F^\leftarrow(p)]\right]\right\}
\end{eqnarray}
and
\begin{equation}\label{auxTh1e1}
\frac{1}{p} g[F^\leftarrow(1-p)] - \frac{1-p}{p}g[F^\leftarrow(p)] \geq g[R(X,p)] \iff p_{R,p}(X) \geq  p_{R,p}[g(X)],
\end{equation}
\begin{equation}\label{auxTh1e2}
\frac{1}{p} g[F^\leftarrow(p)] - \frac{1-p}{p}g[F^\leftarrow(1 - p)] \geq g[L(X,p)] \iff p_{L,p}(X) \leq  p_{L,p}[g(X)].
\end{equation}
\item[f)] If $g(x)$ is continuous and strictly decreasing
\begin{eqnarray}
\label{fL}  p_{L, p}(g(X)) &=& 1 - F \left\{g^\leftarrow\left[\frac{1}{p} g[F^\leftarrow(1 - p)] - \frac{1-p}{p}g[F^\leftarrow(p)]\right]\right\},\\
\label{fR}  p_{R, p}(g(X)) &=& F \left\{g^\leftarrow\left[ \frac{1}{p} g[F^\leftarrow(p)] - \frac{1-p}{p}g[F^\leftarrow(1-p)]\right]\right\},
\end{eqnarray}
and
\begin{equation}\label{auxTh1e3}
\frac{1}{p} g[F^\leftarrow(p)] - \frac{1-p}{p}g[F^\leftarrow(1-p)] \leq g[L(X,p)] \iff  p_{L,p}(X) \leq  p_{R,p}[g(X)],
\end{equation}
\begin{equation}\label{auxTh1e4}
\frac{1}{p} g[F^\leftarrow(1 - p)] - \frac{1-p}{p}g[F^\leftarrow(p)] \leq g[R(X,p)] \iff  p_{L,p}(g(X)) \leq  p_{R,p}(X).
\end{equation}
\item[g)]  $p_{R,p}(X_{(n, n)}) = 1 - F_X^n\left[\frac{1}{p}F^\leftarrow_{X}(\sqrt[n]{1 - p}) - \frac{1-p}{p} F^\leftarrow_{X}(\sqrt[n]{p})\right]$

$p_{L,p}(X_{(n, n)}) = F_X^n\left[\frac{1}{p}F^\leftarrow_{X}(\sqrt[n]{p}) - \frac{1-p}{p} F^\leftarrow_{X}(\sqrt[n]{1 - p})\right]$
\item[h)] $p_{R,p}(X_{(1, n)}) = \left\{1 - F_X\left[\frac{1}{p}F^\leftarrow_{X}(1-\sqrt[n]{p}) - \frac{1-p}{p} F^\leftarrow_{X}(1-\sqrt[n]{1 - p})\right]\right\}^n$

$p_{L,p}(X_{(1, n)}) = 1 - \left\{1 - F_X \left[\frac{1}{p}F^\leftarrow_{X}(1-\sqrt[n]{1 - p}) - \frac{1-p}{p} F^\leftarrow_{X}(1-\sqrt[n]{p})\right]\right\}^n$
\item[i)] For all $t > 0$, $$p_{R,p}(X-t|X>t) = p_{R,p}(X|X>t), \quad p_{L,p}(X-t|X>t) = p_{L,p}(X|X>t).$$
\item[k)] Let $l < u$, and $l, u \in R$. Denote correspondingly the right-, left-,  and double-truncated r.vs. by $X_{RT} = (X|X<u)$, $X_{LT} = (X|X>l)$, $X_{DT} = (X|l<X<u)$. If $0< F(u)$, $F(l) < 1$ and $F(l) \not= F(u)$, then the relations between their probabilities for $p$-outside values, $F^\leftarrow:=F_X^\leftarrow$, and $F:=F_X$ are the following:
    \end{description}
  \begin{eqnarray*}
       p_{L,p}(X_{RT}) &=& \frac{F\left\{\frac{1}{p}F^\leftarrow[pF(u)]-\frac{1-p}{p}F^\leftarrow[(1-p)F(u)]\right\}}{F(u)} \\
       p_{R,p}(X_{RT}) &=& 1-\frac{F\left\{\frac{1}{p}F^\leftarrow[(1-p)F(u)]-\frac{1-p}{p}F^\leftarrow[pF(u)]\right\}}{F(u)}\\
       p_{L,p}(X_{LT}) &=& \frac{F\left\{\frac{1}{p}F^\leftarrow[p+(1-p)F(l)]-\frac{1-p}{p}F^\leftarrow[(1-p)+pF(l)]\right\}-F(l)}{1 - F(l)}\\
       p_{R,p}(X_{LT}) &=& \frac{1-F\left\{\frac{1}{p}F^\leftarrow[1-p+pF(l)]-\frac{1-p}{p}F^\leftarrow[p+(1-p)F(l)]\right\}}{1 - F(l)}\\
       p_{L,p}(X_{DT}) &=& \frac{F\left\{\frac{1}{p}F^\leftarrow[pF(u)+(1-p)F(l)]-\frac{1-p}{p}F^\leftarrow[(1-p)F(u)+pF(l)]\right\}-F(l)}{F(u) - F(l)}\\
       p_{R,p}(X_{DT}) &=& \frac{F(u)-F\left\{\frac{1}{p}F^\leftarrow[(1-p)F(u)+pF(l)]-\frac{1-p}{p}F^\leftarrow[pF(u)+(1-p)F(l)]\right\}}{F(u) - F(l)}.
     \end{eqnarray*}
\item[l)] If $P(X > 0) = 1$, and
\begin{itemize}
  \item $a > 1$, then
\end{itemize}
\begin{eqnarray*}
  p_{L,p}(log_a(X)) &=& F \left\{\sqrt[p]{\frac{F^\leftarrow(p)}{[F^\leftarrow(1 - p)]^{1 - p}}}\right\},\\
  p_{R,p}(log_a(X)) &=& 1 - F\left\{\sqrt[p]{\frac{F^\leftarrow(1 - p)}{[F^\leftarrow(p)^{1 - p}]}}\right\}.
\end{eqnarray*}
\begin{itemize}
  \item $0 < a <  1$, then
\end{itemize}
\begin{eqnarray*}
  p_{L,p}(log_a(X)) &=& 1 - F \left\{\sqrt[p]{\frac{F^\leftarrow(1 - p)}{[F^\leftarrow(p)]^{1 - p}}}\right\},\\
  p_{R,p}(log_a(X)) &=& F\left\{\sqrt[p]{\frac{F^\leftarrow(p)}{[F^\leftarrow(1 - p)]^{1 - p}}}\right\}.
\end{eqnarray*}
\item[m)] For any r.v. $X$, such that $P(X > 0) = 1$, if
\begin{itemize}
  \item If $\alpha <  0$, then
\end{itemize}
\begin{eqnarray*}
  p_{L,p}(X^\alpha) &=& 1 - F\left\{\sqrt[\alpha]{\frac{1}{p}[F^\leftarrow(1 - p)]^\alpha - \frac{1 - p}{p}[F^\leftarrow(p)]^\alpha }\right\},\\
  p_{R,p}(X^\alpha) &=&  F\left\{\sqrt[\alpha]{\frac{1}{p}[F^\leftarrow(p)]^\alpha - \frac{1 - p}{p}[F^\leftarrow(1 - p)]^\alpha }\right\}.
\end{eqnarray*}
\begin{itemize}
  \item If $\alpha >  0$, then 
\end{itemize}
\begin{eqnarray*}
  p_{L,p}(X^\alpha)&=& F\left\{\sqrt[\alpha]{\frac{1}{p}[F^\leftarrow(p)]^\alpha - \frac{1 - p}{p}[F^\leftarrow(1 - p)]^\alpha }\right\},\\
  p_{R,p}(X^\alpha) &=& 1 - F\left\{\sqrt[\alpha]{\frac{1}{p}[F^\leftarrow(1 - p)]^\alpha - \frac{1 - p}{p}[F^\leftarrow(p)]^\alpha }\right\}.
\end{eqnarray*}
\item[n)] If
\begin{itemize}
  \item $a \in (0, 1)$, then
\end{itemize}
\begin{eqnarray*}
  p_{L,p}(a^X) &=& 1 - F\left\{log_a\left[\frac{1}{p}a^{F^\leftarrow(1 - p)} - \frac{1 - p}{p}a^{F^\leftarrow(p)}\right]\right\},\\
  p_{R,p}(a^X) &=&  F\left\{log_a\left[\frac{1}{p}a^{F^\leftarrow(p)} - \frac{1 - p}{p}a^{F^\leftarrow(1 - p)}\right]\right\}.
\end{eqnarray*}
\begin{itemize}
  \item  If $\alpha > 1$, then
\end{itemize}
\begin{eqnarray*}
   p_{L,p}(a^X) &=& F\left\{log_a\left[\frac{1}{p}a^{F^\leftarrow(p)} - \frac{1 - p}{p}a^{F^\leftarrow(1 - p)}\right]\right\},\\
  p_{R,p}(a^X) &=&  1 - F\left\{log_a\left[\frac{1}{p}a^{F^\leftarrow(1 - p)} - \frac{1 - p}{p}a^{F^\leftarrow(p)}\right]\right\}.
\end{eqnarray*}
\end{thm}

\begin{rem} Note that in Theorem 1, l), the expressions for $p_{R,p}(log_a(X))$ and $p_{L,p}(log_a(X))$ do not depend on the exact value of $a$ but only on the fact if  $0 < a <  1$, or $a > 1$. The last means that, according to this classification of absolutely continuous probability distributions with respect to the tails of their c.d.fs., if we decide to change these characteristics and take a logarithm, the exact value of the basis of the logarithm is not important for probabilities for outside values of the transformed distribution. Only the fact that it is bigger or less than $1$ can influence $p_{R,p}(log_a(X))$ and $p_{L,p}(log_a(X))$.
\end{rem}

{\bf Corollary of f):}   Let $p \in (0, 0.5]$ be fixed.
\begin{itemize}
  \item If $P(X > 0) = 1$, then $p_{L,p}(X) \leq  p_{R,p}(\frac{1}{X}).$
  \item If $P(X < 0) = 1$, then $p_{R,p}(X) \leq  p_{L,p}(\frac{1}{X}).$
\end{itemize}

The next corollary corresponds to the well-known experience that taking a logarithm with basis bigger than one of the data we decrease the chance to observe right $p$-outside values and increase the chance to observe left $p$-outside values. Together with Theorem \ref{thm:thm1aux}, they show once again the appropriateness of these characteristics when speaking about the  tail of the observed distribution.

{\bf Corollary of  e):}\label{cor:Cor1e}   Let $p \in (0, 0.5]$ be fixed. Suppose $P(X > 0) = 1$,
\begin{itemize}
  \item if $a >  1$, then
\end{itemize}
\begin{equation}\label{LogarithmsFirstinequalityL}
p_{L,p}(a^X) \leq  p_{L,p}(X) \leq p_{L, p}[log_a(X)],
\end{equation}
\begin{equation}\label{LogarithmsFirstinequalityR}
p_{R,p}[log_a(X)] \leq  p_{R,p}(X) \leq p_{R, p}(a^X).
\end{equation}
\begin{itemize}
  \item If $0 < a <  1$, then
\end{itemize}
\begin{equation}\label{LogarithmsFirstinequalityLalessthan1}
p_{L,p}(a^X) \leq  p_{R, p}(X), \quad p_{L, p}\left(\frac{1}{a^X}\right) \leq p_{L,p}(X) \leq p_{R, p}[log_a(X)],
\end{equation}
\begin{equation}\label{Logarithms1}
p_{L,p}[log_a(X)] \leq  p_{R, p}(X) \leq p_{R, p}\left(\frac{1}{a^X}\right), \quad p_{L,p}(X) \leq  p_{R, p}(a^X).
\end{equation}

According to $p_{R, p}$ characteristics, taking powers bigger than 1 of the data we increase the chance to observe right $p$-outside values, and decrease the chance to observe left $p$-outside values in the observed distribution.

\begin{thm}\label{thm:thm1aux} For $p \in (0, 0.5]$, and $P(X > 0) = 1$,
\begin{description}
\item[a)]  $0 < \alpha_1 \leq \alpha_2$, then
\begin{equation}\label{Cor2ofe1}
p_{L,p}(X^{\alpha_2}) \leq p_{L,p}(X^{\alpha_1}), \quad p_{R,p}(X^{\alpha_1}) \leq p_{R,p}(X^{\alpha_2}).
\end{equation}
\item[b)]  If $\alpha > 1$, and $\frac{1}{1 - p} \leq \left[\frac{F^\leftarrow(1 - p)}{F^\leftarrow(p)}\right]^\alpha$, then $p_{L,p}(X^{\alpha}) = 0$ and
    \begin{equation}\label{Cor2ofe21}
p_{R,p}(X^{1/\alpha}) \leq p_{R,p}(X) \leq p_{R,p}(X^{\alpha}).
\end{equation}
\item[c)] If $\alpha > 1$, and $\frac{1}{1 - p} \geq \left[\frac{F^\leftarrow(1 - p)}{F^\leftarrow(p)}\right]^{\alpha}$, then $p_{R,p}(X^{\alpha}) = 0$ and
  \begin{equation}\label{Cor2ofe2}
p_{L,p}(X^{\alpha}) \leq p_{L,p}(X) \leq p_{L,p}(X^{1/\alpha}).
\end{equation}
\end{description}
\end{thm}

Application of these probabilities requires knowledge about their values for different distributions. Therefore we have calculated some of their explicit forms in the next section.

\section{The most important particular cases}
\label{ParticularCases}

In order to choose the most appropriate class for modeling the tails of the c.d.f. of the observed r.v. we can first calculate the probabilities for left and right $p$-outside values, for as more as possible distributional types, and then to compare these probabilities with corresponding estimators. This approach is analogous to the comparison of the means in cases when we are interested in the center of the distribution.
Let us now present the exact values of these characteristics in some of the most popular cases of probability distributions used in practice for modeling heavy tails. Till the end of this section, we assume that $p \in (0, 0.5]$.

\begin{figure}
\begin{center}
\begin{minipage}[t]{0.45\linewidth}
    \includegraphics[scale=.43]{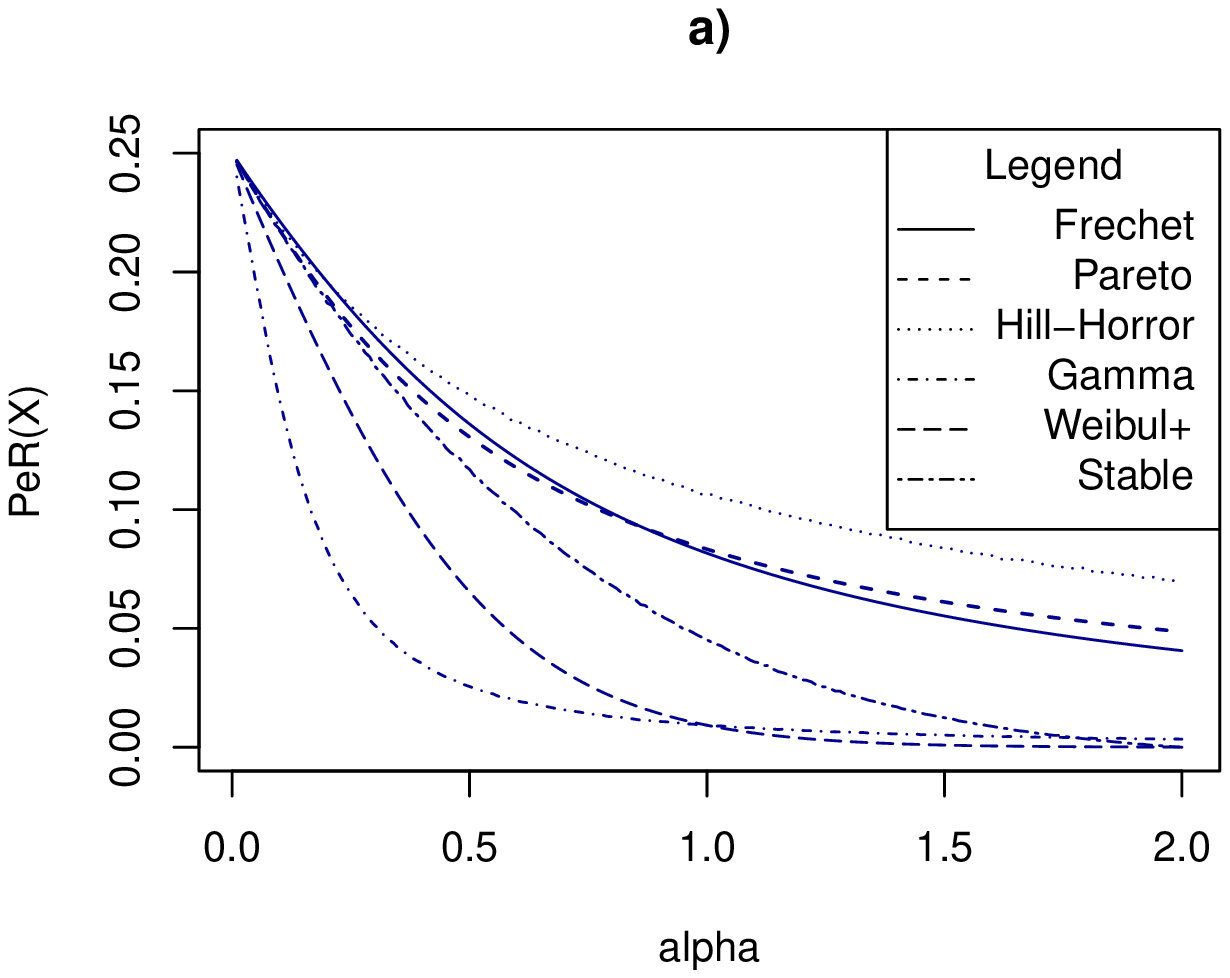}\vspace{-0.3cm}

    \caption{The dependence of $p_{R, 0.25}(X)$ on $\alpha$}
    \label{fig:peRalphaE}
\end{minipage} $ $
\begin{minipage}[t]{0.45\linewidth}
    \includegraphics[scale=.43]{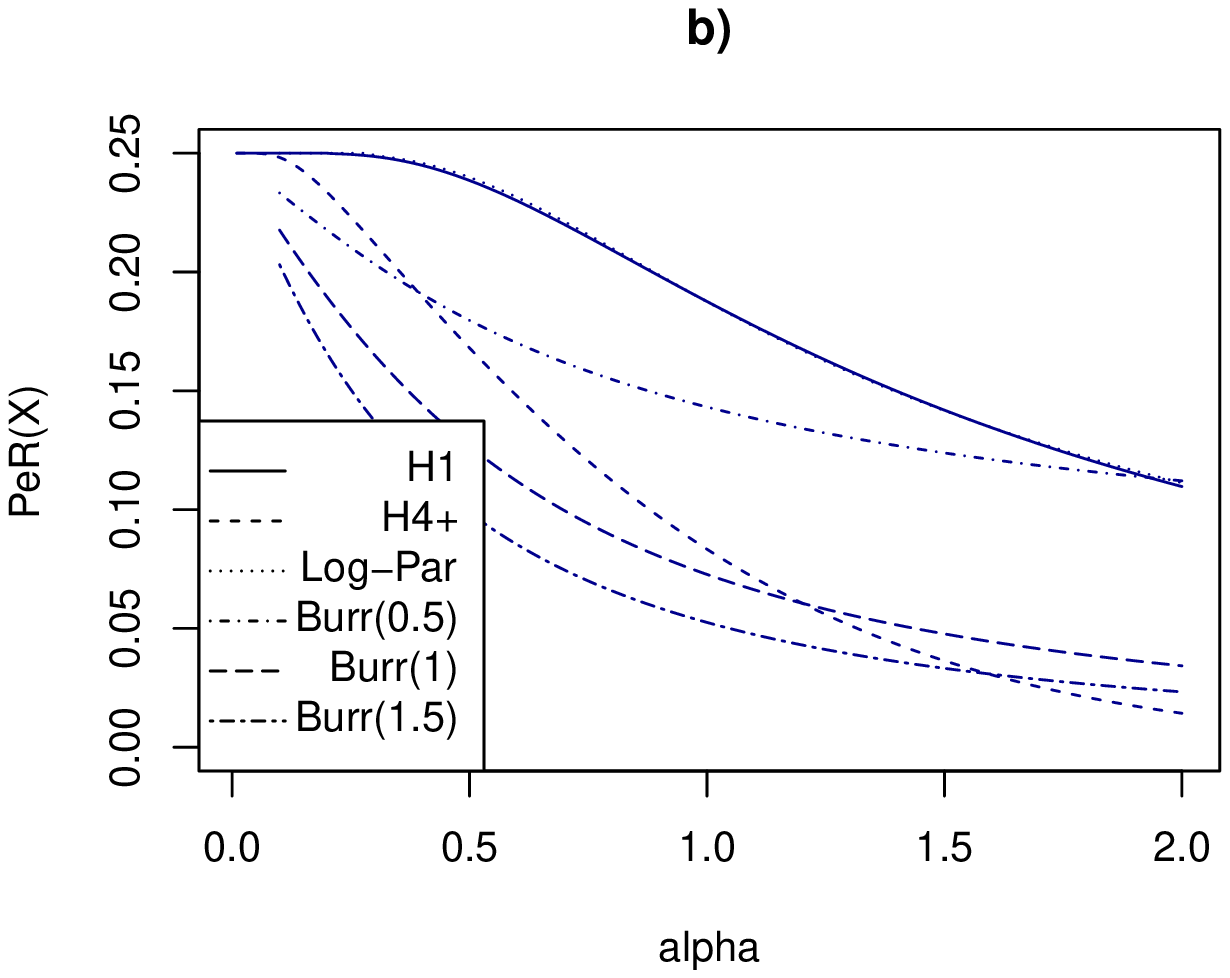}\vspace{-0.3cm}

    \caption{The dependence of $p_{R, 0.25}(X)$ on $\alpha$}
    \label{fig:peRalphaExtremes}
\end{minipage}
\end{center}
\end{figure}

The dependencies of $p_{R, 0.25}$ on the parameter $\alpha$, which characterises the tail of the corresponding distribution in cases when $F$ is $Gamma(\alpha, \beta)$, or Fr$\acute{e}$chet, Pareto, Stable, Weibull positive, $H_1$, $H_4$, log-Pareto, Hill horror or Burr distributed are depicted on  Figure \ref{fig:peRalphaE} and Figure \ref{fig:peRalphaExtremes}, and could be seen also in Jordanova and Petkova (2018) \cite{MoniPoli2018}, and in the supplementary material of Soza et al. (2019) \cite{JordanovaStehlik2018}.

\begin{itemize}
\item Exponential distribution. Let $\lambda > 0$, and $X$ be Exponential with mean $\frac{1}{\lambda}$.
\end{itemize}
 It is well known that $\lambda$ is a scale parameter of the exponential distribution, therefore due to Th. 1, c) without lost of generality (w.l.g.) we can assume that $\lambda = 1$ and this will not change the values of $p_{L, p}(X)$ and $p_{R, p}(X)$. In this case,
 $$ p_{L,p}(X) = P\left[X < log\frac{p^{\frac{1-p}{p}}}{(1-p)^{\frac{1}{p}}}\right] = \left\{
 \begin{array}{ccc}
   0 & , & p \in (0, p_0), \\
   p\frac{(1-p)^{\frac{1}{p}}}{p^{\frac{1}{p}}} & , & p \in [p_0, \frac{1}{2})
 \end{array}
 \right.,$$
   where $p_0$ is the solution of the equation $(1-p_0)log(p_0)=log(1-p_0)$ and $p_0 \approx 0.4096$.
 $$p_{R, p}(X) = P\left\{X > log\left[\left(\frac{1 - p}{p}\right)^{\frac{1}{p}}\frac{1}{1-p}\right]\right\} = (1-p)\left(\frac{p}{1-p}\right)^{\frac{1}{p}}.$$

In particular, for  $p_{L, 0.25}(X) = 0$, the empirical right fence is asymptotically unbiased and efficient estimator for the theoretical right fence
$$\lim_{i \to \infty} ER_{4i-1}(X,0.25) = R(X,0.25) = 2log(2) + 3log(3) = log(108)$$
$$\lim_{i \to \infty} DR_{4i-1}(X,0.25) = \lim_{i \to \infty} [16\psi'(i) - \psi'(4i) - 15\psi'(3i)] = 0.$$
where $\psi'(i) = \frac{\partial^2 log \Gamma(z)}{\partial z^2}$ is the Polygamma function (for the last limit see Guo and Feng (2013) \cite{GuoFeng2013}), and $p_{R, 0.25}(X) = \frac{1}{108} = 0.00925(925).$
See Jordanova and Petkova (2017-2018) \cite{MoniPoli2018,MoniPoli2017}.

Further on in this section (due to properties b) and c) Theorem \ref{thm:thm1}, w.l.g. we assume that $\mu = 0$ and $\sigma = 1$.

\begin{itemize}
\item Generalized Pareto distribution (GPD). Consider $\mu \in R$ and $\sigma > 0$.
              \end{itemize}
\begin{eqnarray}\label{GPD}
F_X(x) = \left\{\begin{array}{ccc}
                     1 - (1 + \xi \frac{x - \mu}{\sigma})^{-1/\xi} &, & x > \mu, \xi > 0 \\
                     1 - (1 + \xi \frac{x - \mu}{\sigma})^{-1/\xi} &, & \mu \leq x \leq \mu - \frac{\sigma}{\xi}, \xi < 0 \\
                     1- e^{-\frac{x - \mu}{\sigma}} & , & x > \mu, \xi \to 0.
                   \end{array}
 \right.
 \end{eqnarray}
 We have already considered the case $\xi \to 0$. In that case, it is well known that, the GPD coincides with Exponential distribution. So, here we assume that $\xi \not = 0$. Then the quantile function is
$F_X^\leftarrow(p) = \frac{1}{\xi}[(1-p)^{-\xi}-1].$ We replace it in the formula for $L(X; p)$, then in the definition for $p_{L, p}(X)$, and  obtain
  $$ p_{L,p}(X) = P\left\{X < \frac{1}{\xi}\left[\frac{1}{p}(1-p)^{-\xi} -1 + p^{-\xi} - p^{-\xi-1}\right]\right\}.$$
 In order to replace (\ref{GPD}) in the last probability we need to consider separately the following two cases:
 \begin{itemize}
   \item[$\cdot$] Case $\xi > 0$. In this case
   $\frac{1}{\xi}\left[\frac{1}{p}(1-p)^{-\xi} -1 + p^{-\xi} - p^{-\xi-1}\right] < \frac{1}{\xi}$ $\left[\frac{1}{p}-1 +\right.$ $\left. p^{-\xi} - p^{-\xi-1} \right]$. The last expression is equal to  $ \frac{1}{\xi} \left(1- \frac{1}{p}\right)\left(\frac{1}{p^\xi}-1\right) < 0 $, therefore $p_{L,p}(X) = 0$.
   \item[$\cdot$] Case $\xi < 0$. In this case $\frac{1}{\xi}\left[\frac{1}{p}(1-p)^{-\xi} -1 + p^{-\xi} - p^{-\xi-1}\right] > 0$, therefore
 \end{itemize}
    \begin{eqnarray*}
  p_{L,p}(X) &=& P\left\{X < \frac{1}{\xi}\left[\frac{1}{p}(1-p)^{-\xi} -1 + p^{-\xi} - p^{-\xi-1}\right]\right\}\\
  &=& 1 - (\frac{1}{p}(1-p)^{-\xi} + p^{-\xi} - p^{-\xi-1})^{-1/\xi}
 \end{eqnarray*}

 Analogously we replace the quantile function in the definition for $R(X; p)$, then in $p_{R, p}(X)$, and  obtain
 \begin{eqnarray*}
  p_{R, p}(X) &=& P\left\{X > \frac{1}{\xi}[p^{-\xi}-1] + \frac{1-p}{p}\left[\frac{1}{\xi}(p^{-\xi}-1) - \frac{1}{\xi}[(1-p)^{-\xi}-1]\right]\right\}\\
 &=& P\left\{X > \frac{1}{\xi}\left[\frac{1-p}{p}[1-(1-p)^{-\xi}]-\frac{1}{p}\left(1 - p^{-\xi}\right)\right]\right\}
 \end{eqnarray*}

In order to calculate this expression we need to determine the sign of $\frac{1}{\xi}\left[\frac{1-p}{p}[1-(1-p)^{-\xi}]-\frac{1}{p}\left(1 - p^{-\xi}\right)\right]$. Therefore again we consider two cases.
 \begin{itemize}
   \item[$\cdot$] Case $\xi > 0$. Because of $p \in (0, 0.5]$, we have $p^{-\xi} > (1-p)^{-\xi}$,  $p(1-p)^{-\xi} > p$ and
\end{itemize}
$$(1-p)^{-\xi} + p  < p^{-\xi} + p(1-p)^{-\xi}$$
$$1-(1-p)^{-\xi}-p[1-(1-p)^{-\xi}] > 1 - p^{-\xi}$$
$$(1-p)[1-(1-p)^{-\xi}] > 1 - p^{-\xi}$$
$$\frac{1-p}{p}[1-(1-p)^{-\xi}] > \frac{1}{p}\left(1 - p^{-\xi}\right)$$
Therefore
\begin{eqnarray*}
  p_{R, p}(X) &=& P\left\{X > \frac{1}{\xi}\left[\frac{1-p}{p}[1-(1-p)^{-\xi}]-\frac{1}{p}\left(1 - p^{-\xi}\right)\right]\right\}\\
 &=& \left\{1 + \frac{1-p}{p}[1-(1-p)^{-\xi}]-\frac{1}{p}\left(1 - p^{-\xi}\right)\right\}^{-1/\xi}\\
 &=& \left\{\frac{p^{-\xi}}{p}+\left(1 - \frac{1}{p}\right)(1 - p)^{-\xi}\right\}^{-1/\xi}
 \end{eqnarray*}

In case $p = 0.25$, and $\alpha = \frac{1}{\xi}$ this expression coincide with the one in Jordanova and Petkova (2018)\cite{MoniPoli2018}
\begin{equation}\label{PR025Pareto}
p_{R, 0.25}(X) = \frac{3}{4(4.3^{\frac{1}{\alpha}}-3)^{\alpha}}.
\end{equation}

\begin{itemize}
   \item[$\cdot$] Case $\xi < 0$. Because of $p \in (0, 0.5)$  we have that $\frac{1-p}{p} < \frac{1}{p}$, $1-(1-p)^{-\xi} < 1 - p^{-\xi}$  and $\frac{1}{\xi}\left[\frac{1-p}{p}[1-(1-p)^{-\xi}]-\frac{1}{p}\left(1 - p^{-\xi}\right)\right] > 0.$
\end{itemize}

  In case
 $$\frac{1}{\xi}\left[\frac{1-p}{p}[1-(1-p)^{-\xi}]-\frac{1}{p}\left(1 - p^{-\xi}\right)\right] > -\frac{1}{\xi}$$
   $$(1-p)[1-(1-p)^{-\xi}] - 1 + p^{-\xi} < -p$$
  $$p^{-\xi} < (1 - p)^{1-\xi}$$
    $$\xi  < \frac{log(1 - p)}{log\left(\frac{1 - p}{p}\right)}$$
  we have that $p_{R, p}(X) =0$.

  If $0 > \xi  > \frac{log(1 - p)}{log\left(\frac{1 - p}{p}\right)}$,
\begin{eqnarray*}
p_{R, p}(X) &=& P\left\{X > \frac{1}{\xi}\left[\frac{1-p}{p}[1-(1-p)^{-\xi}]-\frac{1}{p}\left(1 - p^{-\xi}\right)\right]\right\}\\
&=& \left\{1 + \frac{1-p}{p}[1-(1-p)^{-\xi}]-\frac{1}{p}(1 - p^{-\xi})\right\}^{-\frac{1}{\xi}} \\
&=&  \left\{(1-p)^{-\xi}\left(1 - \frac{1}{p}\right) + \frac{1}{p} p^{-\xi}\right\}^{-\frac{1}{\xi}}.
\end{eqnarray*}

When speaking about heavy tails we can not forget about Extreme value distributions with respect to linear transformations. Therefore in the next three points, we will consider them. At the beginning of the last century Fisher and Tippet (1928)\cite{FisherTippet}, Gnedenko (1943) \cite{Gnedenko1943}, and  Gumbel (1958)\cite{Gumbel1958} have shown that they appear as limiting distributions of maxima of i.i.d. r.vs. after appropriate affine transformations.
\begin{itemize}
\item  Fr$\acute{e}$chet distribution.  Let $\alpha > 0$. W.l.g. we assume that in (\ref{Frechet}) $c = 1$.
 \end{itemize}
 Therefore for $p \in (0, 1)$, $F^\leftarrow(p) = (-log(p))^{-\frac{1}{\alpha}}$.
 $$ p_{L,p}(X) =  P\left\{X < (-log(p))^{-\frac{1}{\alpha}} - \frac{1-p}{p}[(-log(1 - p))^{-\frac{1}{\alpha}} - (-log(p))^{-\frac{1}{\alpha}}]\right\}$$
 and because of $(-log(p))^{-\frac{1}{\alpha}} > \frac{1-p}{p}[(-log(1 - p))^{-\frac{1}{\alpha}} - (-log(p))^{-\frac{1}{\alpha}}]$ $\Longleftrightarrow $
 $$-log(p) < (1-p)^{-\alpha}(-log(1-p)) \Longleftrightarrow \alpha_0: = -\frac{log\left(\frac{log(p)}{log(1-p)}\right)}{log(1-p)} < \alpha$$
 $$p_{L, p}(X) = exp\left\{-\left[\frac{1}{p}(-log\,p)^{-1/\alpha}+ \left(1-\frac{1}{p}\right)[-log\,(1-p)]^{-1/\alpha}\right]^{-\alpha}\right\},$$ when $\alpha > \alpha_0$, and $p_{L, p}(X) =0$, for $\alpha \in (0, \alpha_0]$.

 Analogously
\begin{eqnarray*}
  p_{R, p}(X) &=&  P\left\{X > (-log(1 - p))^{-\frac{1}{\alpha}} + \frac{1-p}{p}[(-log(1 - p))^{-\frac{1}{\alpha}} - (-log(p))^{-\frac{1}{\alpha}}]\right\}\\
&=&  P\left[X >  \frac{1}{p}(-log(1 - p))^{-\frac{1}{\alpha}} - \frac{1-p}{p} (-log(p))^{-\frac{1}{\alpha}}\right].
\end{eqnarray*}

Now we need to consider the expression after the inequality. As far as for all $p \in (0, 0.5)$ the following four expressions are equivalent
$$(-log(1 - p))^{-\frac{1}{\alpha}} > (1-p)(-log(p))^{-\frac{1}{\alpha}}$$
$$\frac{-log(1 - p)}{-log(p)} < (1-p)^{-\alpha}$$
$$log\left[\frac{-log(1 - p)}{-log(p)}\right] < -\alpha \,\,log(1-p)$$
$$-\frac{log\left[\frac{-log(1 - p)}{-log(p)}\right]}{log(1-p)} < 0 < \alpha $$
we have
$$p_{R, p}(X) = exp\left\{-\left(\frac{1}{p}(-log(1 - p))^{-\frac{1}{\alpha}} - \frac{1-p}{p} (-log(p))^{-\frac{1}{\alpha}}\right)^{-\alpha}\right\}.$$

Figure \ref{fig:peRalphaE} represents the dependence of $p_{R, p}(X)$ on $\alpha$ in case $p = 0.25$. The explicit formula for this case can be seen also in Jordanova and Petkova (2018) \cite{MoniPoli2018}.

\begin{itemize}
\item  Weibull negative distribution.  Consider $\alpha > 0$. W.l.g. $X \in NegWeibull(\alpha,$ $1,$ $0)$. See (\ref{Weibull}). Therefore for $p \in (0, 1)$, $F^\leftarrow(p) = -(-log(p))^{\frac{1}{\alpha}}.$
\end{itemize}
Note that its positive version $-X$ coincides in distribution with a standard Exponentially distributed r.v. raised to the power $\frac{1}{\alpha}$. Therefore, for $\alpha > 1$,  according to our classification, this distribution has heavier right tail than the exponential one and for $\alpha < 1$ vice versa. This can be seen also on Figure \ref{fig:peRalphaE}, where standard exponential distribution is depicted as $\Gamma(1,1)$. More precisely
 \begin{eqnarray*}
  p_{L,p}(X) &=& P\left\{X < -[-log(p)]^{\frac{1}{\alpha}} - \frac{1-p}{p}\left[[-log(p)]^{\frac{1}{\alpha}} - [-log(1 - p)]^{\frac{1}{\alpha}}\right]\right\}\\
  &=& P\left\{X < \frac{1-p}{p}[-log(1 - p)]^{\frac{1}{\alpha}} - \frac{1}{p}[-log(p)]^{\frac{1}{\alpha}} \right\}
 \end{eqnarray*}
 and because of for all $p \in (0, 0.5]$ we have $(1-p)[-log(1 - p)]^{\frac{1}{\alpha}} < [-log(p)]^{\frac{1}{\alpha}}$, therefore
 \begin{equation}\label{WNalpha}
 p_{L, p}(X) = exp\left\{-\left[\frac{1}{p}[-log(p)]^{\frac{1}{\alpha}} - \frac{1-p}{p}[-log(1 - p)]^{\frac{1}{\alpha}}\right]^{   \alpha}\right\}
 \end{equation}

Figure \ref{fig:peRalphaE}, depicts the dependence of $p_{L, 0.25}(X) = p_{R, 0.25}(-X)$, (i.e. $-X$ is Weibull positive)  on $\alpha$.

 Analogously
\begin{eqnarray*}
  p_{R, p}(X) &=& P\left\{X > -(-log(1 - p))^{\frac{1}{\alpha}} + \frac{1-p}{p}[(-log(p))^{\frac{1}{\alpha}} -(-log(1 - p))^{\frac{1}{\alpha}} ]\right\}\\
&=&  P\left[X >  \frac{1-p}{p}(-log(p))^{\frac{1}{\alpha}} - \frac{1}{p} (-log(1 - p))^{\frac{1}{\alpha}}\right]
\end{eqnarray*}

As far as for all $p \in (0, 0.5]$
$$(-log(1 - p))^{\frac{1}{\alpha}} < (1-p)(-log(p))^{\frac{1}{\alpha}}$$
$$\frac{-log(1 - p)}{-log(p)} > (1-p)^{-\alpha}$$
$$log\left[\frac{-log(1 - p)}{-log(p)}\right] > \alpha \,\,log(1-p)$$
$$\alpha_1: = \frac{log\left[\frac{-log(1 - p)}{-log(p)}\right]}{log(1-p)} > \alpha $$
we have
$$p_{R, p}(X) = \left\{\begin{array}{ccc}
                        0 & , & \alpha \in (0, \alpha_1] \\
                        1 - exp\left\{-\left(\frac{1}{p}(-log(1 - p))^{\frac{1}{\alpha}} - \frac{1-p}{p} (-log(p))^{\frac{1}{\alpha}}\right)^\alpha\right\}& , & \alpha > \alpha_1                      \end{array}
 \right..$$

\begin{itemize}
\item  Gumbell distribution.  Let $\alpha > 0$, $\mu \in R$ and $\sigma > 0$
\end{itemize}
$$F_X(x) = exp\left\{-exp\left[-\frac{x-\mu}{\sigma}\right]\right\},\quad  x \in R.$$
W.l.g. $\mu = 0$ and $\sigma = 1$. $F^\leftarrow(p) = -log(-log(p))$,
\begin{eqnarray*}
  p_{L,p}(X) &=& P\left\{X < -log(-log(p)) - \frac{1-p}{p}[log(-log(p)) - log(-log(1 - p))]\right\}\\
  &=& P\left\{X < \frac{1-p}{p}log(-log(1 - p)) - \frac{1}{p}log(-log(p))\right\}  = p^{\left[\frac{log(p)}{log\,(1 - p)}\right]^\frac{1}{p}}
 \end{eqnarray*}
\begin{eqnarray*}
  p_{R, p}(X) &=&  P\left\{X > -log(-log(1 - p)) + \frac{1-p}{p}[log(-log(p))-log(-log(1 - p))]\right\}\\
  &=& P\left\{X > \frac{1-p}{p}log(-log(p)) - \frac{1}{p}log(-log(1 - p))\right\}  = 1 - p^{\left[\frac{log(1-p)}{log\,(p)}\right]^\frac{1}{p}}
\end{eqnarray*}

Jordanova and Petkova (2018) \cite{MoniPoli2018} have calculated that $p_{L, 0.25}(X)  \approx 4.264\times10^{-68}$ and $p_{R, 0.25}(X) \approx 0.002568$.

\begin{itemize}
\item Logistic distribution. Assume $\mu \in R$, $\sigma > 0$ and
\end{itemize}
$$F_X(x) = \frac{1}{1+exp\left(-\frac{x-\mu}{\sigma}\right)},\quad  x \in R.$$
W.l.g. $\mu = 0$ and $\sigma = 1$. $F^\leftarrow(p) = log\left(\frac{p}{1 - p}\right)$,
$$p_{R, p}(X) = p_{L,p}(X) = P\left(X < \frac{2-p}{p}log\frac{p}{1 - p}\right) = \left[1 + \left(\frac{1-p}{p}\right)^{\frac{2-p}{p}}\right]^{-1}$$
 $\leq \left[1 + 3^5\right]^{-1} = p_{R, 0.25}(X)\approx 0.000457,  p \in (0, 0.25).$

\begin{figure}
\begin{center}
\begin{minipage}[t]{0.45\linewidth}
    \includegraphics[scale=.43]{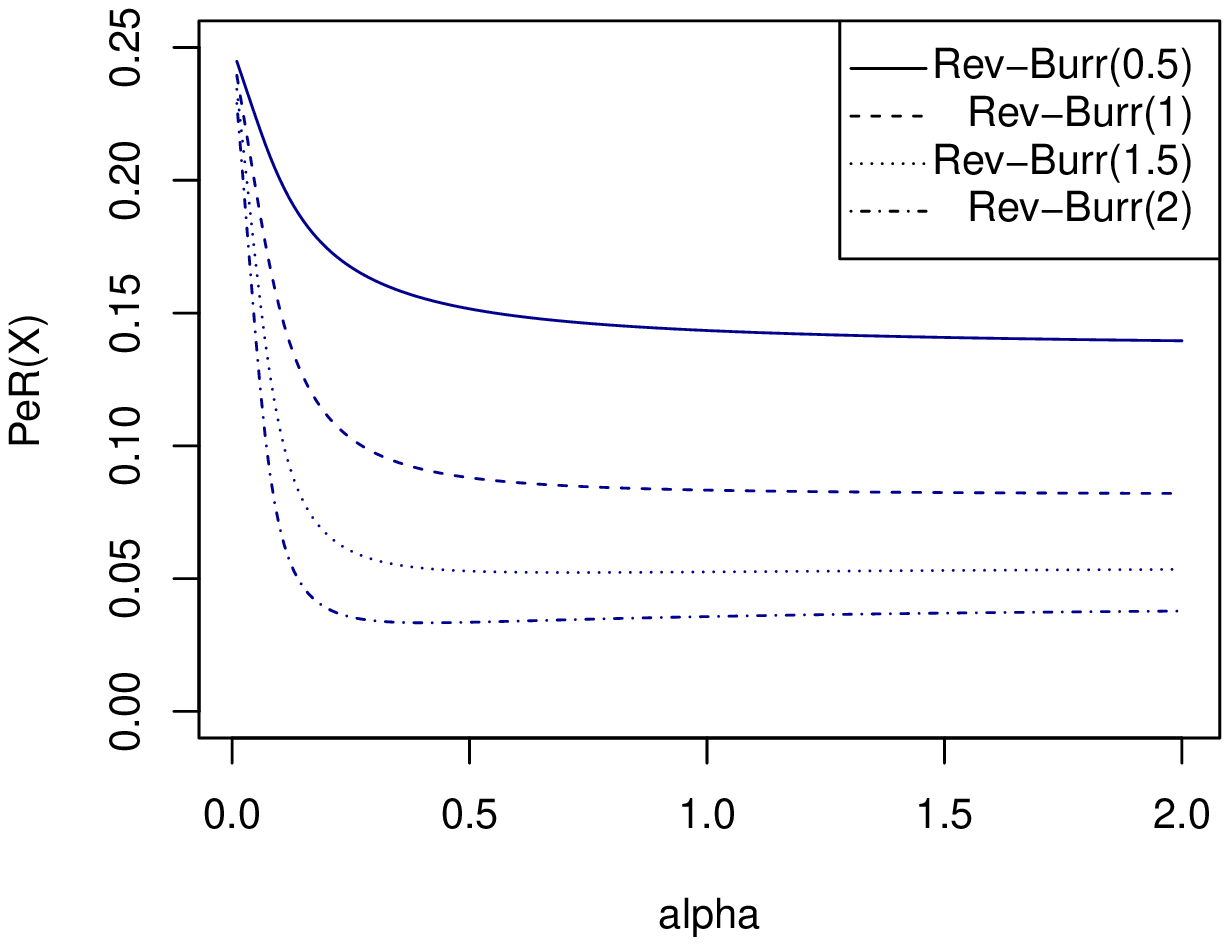}\vspace{-0.3cm}

    \caption{The dependence of $p_{R, 0.25}(-X)$ on parameter $\alpha$ in Reverse-Burr case for $ $ $\tau = 0.5, 1, 1.5, 2$}
    \label{fig:peRalphaExtremes3}
\end{minipage} $ $
\begin{minipage}[t]{0.45\linewidth}
    \includegraphics[scale=.43]{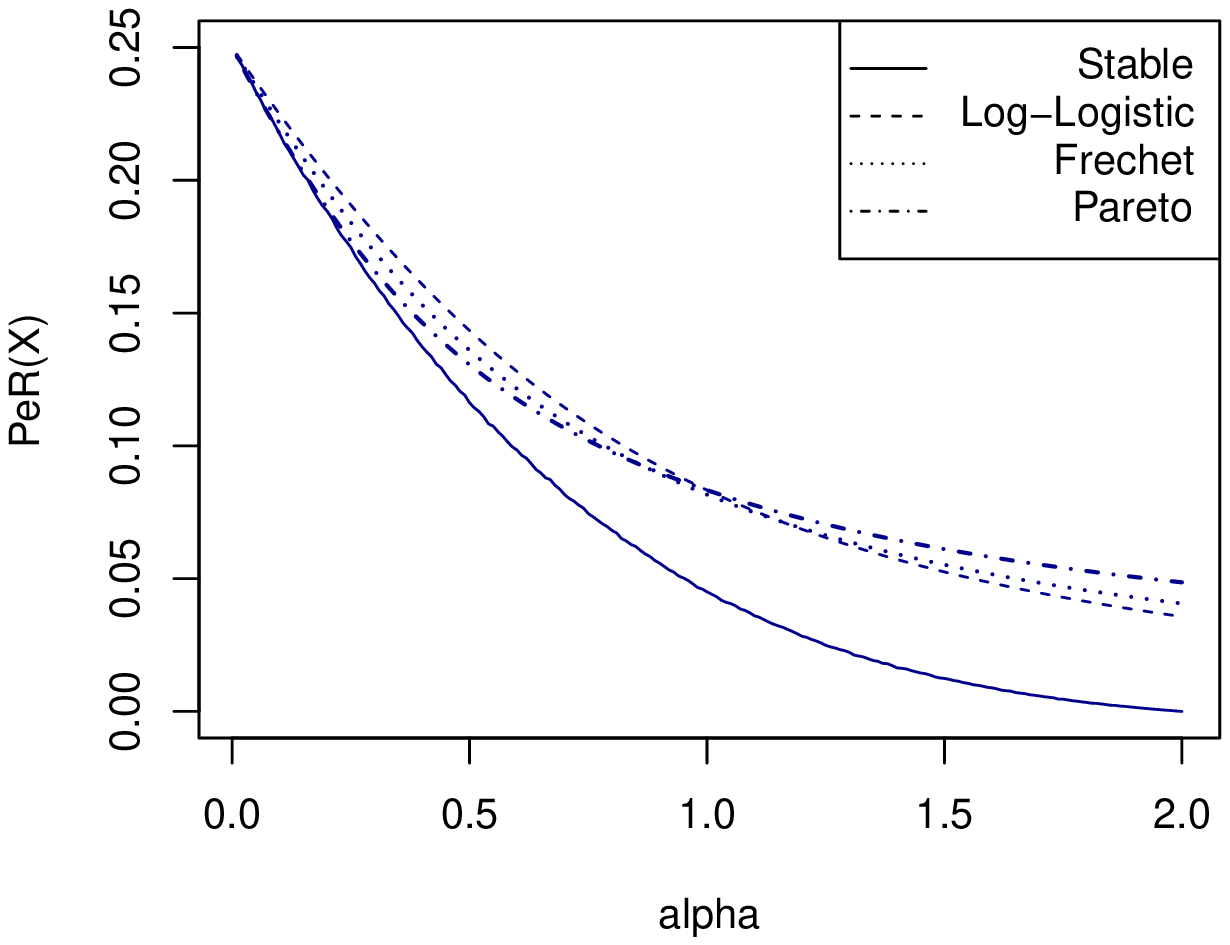}\vspace{-0.3cm}

    \caption{The dependence of $p_{R, 0.25}(X)$ on $\alpha$}
    \label{fig:peRalphaLogLogistic}
\end{minipage}
\end{center}
\end{figure}

\begin{itemize}
\item Log-logistic distribution. Assume $\mu \in R$, $\sigma > 0$ and
\end{itemize}
\begin{equation}\label{LolLogidticCDF}
F_X(x) = \frac{1}{1+\left(\frac{x-\mu}{\sigma}\right)^{-\alpha}},\quad  x > \mu.
\end{equation}
W.l.g. $\mu = 0$ and $\sigma = 1$. As far as $F^\leftarrow(p) = \sqrt[\alpha]{\frac{p}{1 - p}}$, and for $p \in (0, 0.5]$
$$R(X,p) = \frac{1}{p}\sqrt[\alpha]{\frac{1-p}{p}}-\frac{1 - p}{p}\sqrt[\alpha]{\frac{p}{1 -p}} > 0$$
therefore by definition
\begin{equation}\label{pRpLogLogistic}
 p_{R, p}(X) = \frac{1}{1+\left(\frac{1}{p}\sqrt[\alpha]{\frac{1-p}{p}}-\frac{1 - p}{p}\sqrt[\alpha]{\frac{p}{1 -p}}\right)^{\alpha}}.
\end{equation}
The plot of this function of $\alpha > 0$ for $p = 0.25$ can be seen on Figure \ref{fig:peRalphaLogLogistic}. We observe that for a fixed $\alpha > 0$ its right tail is very similar to the corresponding tails of Pareto and Fr$\acute{e}$chet distribution and heavier than the Stable one.
Analogously
$$L(X,p) = \frac{1}{p}\sqrt[\alpha]{\frac{p}{1-p}}-\frac{1 - p}{p}\sqrt[\alpha]{\frac{1-p}{p}},$$
 therefore for $0 < \alpha < 2 [log_{1-p}(p) - 1]$, and $p \in (0, 0.5]$
$$p_{L, p}(X) = 1 - \frac{1}{1+\left(\frac{1}{p}\sqrt[\alpha]{\frac{p}{1 - p}}-\frac{1 - p}{p}\sqrt[\alpha]{\frac{1 - p}{p}}\right)^{\alpha}}.$$
and $p_{L,p}(X) = 0$ otherwise.

In the next two cases we assume that $\tau > 0$ and $\alpha > 0$. W.l.g. $\mu = 0$, because it is a location parameter, and $\sigma = 1$ and $\delta = 1$ because they are scale parameters. See Burr (1942) \cite{Burr} or Einmahl et al. (2008) \cite{EinmahlGuillou}.

\begin{itemize}
\item Burr distribution.  Let
\end{itemize}
\begin{equation}\label{f}
F_X(x) =\left\{\begin{array}{ccc}
                   0 & ,  & x < \mu \\
                   1 - \left[\frac{\delta}{\delta + \left(\frac{x-\mu}{\sigma}\right)^{\tau}}\right]^\alpha & , & x \geq \mu
                 \end{array}
\right..
\end{equation}
For $p \in (0, 1)$, the quantile function is
$$F^\leftarrow(p) = \sqrt[\tau]{(1 - p)^{-1/\alpha}-1} $$
(see also Nair et al. (2013) \cite{Nair2013}). Therefore in our context from (\ref{f}) we have that if
$$\frac{(1-p)^{-\frac{1}{\alpha}}-1}{p^{-\frac{1}{\alpha}}-1} > (1-p)^\tau,$$
then
\begin{eqnarray*}
p_{L,p}(X) &=&  P\left\{X < \frac{1}{p} \sqrt[\tau]{(1 - p)^{-1/\alpha}-1} - \frac{1-p}{p}\sqrt[\tau]{p^{-1/\alpha}-1}\right\}\\
&=&  1 - \left\{\frac{p^\tau}{p^\tau + \left[\sqrt[\tau]{(1 - p)^{-1/\alpha}-1} - (1-p)\sqrt[\tau]{p^{-1/\alpha}-1}\right]^{\tau}}\right\}^\alpha
\end{eqnarray*}
and $p_{L,p}(X) = 0$ otherwise.

For $p \in (0, 0.25)$ and $\tau > 0$ the inequality $$ \frac{p^{-\frac{1}{\alpha}}-1}{(1-p)^{-\frac{1}{\alpha}}-1} \geq 1 \geq (1-p)^\tau,$$ the definition of $p_{R, p}(X)$, and (\ref{f}) entail
\begin{eqnarray*}
 p_{R, p}(X) &=&  P\left[X > \frac{1}{p} \sqrt[\tau]{p^{-1/\alpha}-1} - \frac{1-p}{p}\sqrt[\tau]{(1 - p)^{-1/\alpha}-1}\right]\\
&=&  \left[\frac{p^\tau}{p^\tau + \left(\sqrt[\tau]{p^{-1/\alpha}-1} - (1-p)\sqrt[\tau]{(1 - p)^{-1/\alpha}-1}\right)^{\tau}}\right]^\alpha.
\end{eqnarray*}

The dependence of $p_{R,0.25}(X)$ on $\alpha$ and for different values of $\tau$ is depicted on Figure \ref{fig:peRalphaExtremes}. We see that when $\tau$ increases, the chance to observe $p$-outside values in the considered distribution decreases. The last means that for Burr distribution not only $\alpha$ but also $\tau$ influences the tail-behaviour.

\begin{itemize}
\item Reverse Burr distribution.
\end{itemize}
\begin{equation}\label{f1}
F_X(x) =\left\{\begin{array}{ccc}
                   1 - \left[\frac{\delta}{\delta + \left(\frac{\mu-x}{\sigma}\right)^{-\tau}}\right]^\alpha & , & x \leq \mu\\
                   1 & ,  & x > \mu
                 \end{array}
\right..
\end{equation}

For $p \in (0, 1)$, the corresponding quantile function is
$$F^\leftarrow(p) = -[(1 - p)^{-1/\alpha}-1]^{-1/\tau}.$$

We are interested in the case when $p \in (0, 0.5]$. It guarantees that
$$(1-p)[p^{-1/\alpha}-1]^{-1/\tau} < [(1 - p)^{-1/\alpha}-1]^{-1/\tau}.$$

Therefore from (\ref{f1}) we have
\begin{eqnarray*}
p_{L,p}(X) &=&  P\left\{X < \frac{1-p}{p}[p^{-1/\alpha}-1]^{-1/\tau} - \frac{1}{p} [(1 - p)^{-1/\alpha}-1]^{-1/\tau}\right\}\\
&=& 1 - \left\{\frac{1}{1 + p^\tau\left[[(1 - p)^{-1/\alpha}-1]^{-1/\tau}-(1-p)(p^{-1/\alpha}-1)^{-1/\tau}\right]^{-\tau}}\right\}^\alpha.
\end{eqnarray*}

The dependence of $p_{R,0.25}(-X) = p_{L,0.25}(X)$ on $\alpha$ and for different values of $\tau$ is depicted on Figure \ref{fig:peRalphaExtremes3}.
We observe that for $\alpha > 0.5$ the tail behavior of the Reverse-Burr distribution is much more sensitive on $\tau$ than on $\alpha$.

If $(1-p)[(1 - p)^{-1/\alpha}-1]^{-1/\tau} < [p^{-1/\alpha}-1]^{-1/\tau},$
\begin{eqnarray*}
 p_{R, p}(X) &=&  P\left[X > \frac{1-p}{p}[(1 - p)^{-1/\alpha}-1]^{-1/\tau} - \frac{1}{p} [p^{-1/\alpha}-1]^{-1/\tau} \right]\\
&=& \left[\frac{1}{1 + p^\tau\left((1-p)[(1 - p)^{-1/\alpha}-1]^{-1/\tau} - [p^{-1/\alpha}-1]^{-1/\tau}\right)^{-\tau}}\right]^\alpha.
\end{eqnarray*}
and $p_{R, p}(X) = 0$ otherwise.

\begin{itemize}
\item Gompertz distribution. Assume $\mu \in R$, $\sigma > 0$, $\alpha > 0$, and
\end{itemize}
\begin{equation}\label{f11}
F_X(x) =\left\{\begin{array}{ccc}
                    0 & , & x \leq \mu\\
                   1 - exp\left\{-\alpha\left[exp\left(\frac{x-\mu}{\sigma}\right)-1\right]\right\} & ,  & x > \mu
                 \end{array}
\right..
\end{equation}

For $\mu = 0$, $\sigma = 1$ and $p \in (0, 1)$,
$F^\leftarrow(p) = log\left[1 - \frac{log(1 - p)}{\alpha}\right].$

For $1-\frac{1}{\alpha}log(1 - p) \geq \left[1-\frac{1}{\alpha}log(p)\right]^{1-p}$ we have that $p_{L,p}(X)$ is equal to
 \begin{eqnarray*}
   & & P\left\{X < \frac{1}{p}log\left[1 - \frac{log(1 - p)}{\alpha}\right] - \frac{1 - p}{p}log\left[1 - \frac{log(p)}{\alpha}\right]\right\}\\
  &=& 1 - exp\left\{-\alpha\left[exp\left(\frac{1}{p}log\left[1 - \frac{log(1 - p)}{\alpha}\right] - \frac{1 - p}{p}log\left[1 - \frac{log(p)}{\alpha}\right]\right)-1\right]\right\}\\
   &=& 1 - exp\left\{-\alpha\left[\sqrt[p]{\frac{1-\alpha^{-1}log(1 - p)}{[1-\alpha^{-1}log(p)]^{1-p}}}-1\right]\right\}
 \end{eqnarray*}
 and $p_{L, p}(X) = 0$ otherwise. In particular $p_{L, 0.25}(X) = 0$.

 As far as for all $p \in (0, 0.5]$
 $\frac{log[1-\frac{1}{\alpha}log(p)]}{log[1-\frac{1}{\alpha}log(1 - p)]} > 1 - p,$ therefore from (\ref{f11}) and the definition of $p_{R, p}(X)$ we have that it is equal to
\begin{eqnarray*}
   & & P\left\{X > \frac{1}{p}log\left[1 - \frac{log(p)}{\alpha}\right]-\frac{1-p}{p}log\left[1 - \frac{log(1 - p)}{\alpha}\right]\right\}\\
&=&  exp\left\{-\alpha\left[exp\left[\frac{1}{p}log\left[1 - \frac{log(p)}{\alpha}\right]-\frac{1-p}{p}log\left[1 - \frac{log(1 - p)}{\alpha}\right]\right]-1\right]\right\}\\
&=& exp\left\{-\alpha\left[\sqrt[p]{\frac{1-\alpha^{-1}log(p)}{[1-\alpha^{-1}log(1 - p)]^{1-p}}}-1\right]\right\} \approx 0.
\end{eqnarray*}

Note that or all fixed $\alpha > 0$, and within the considered distributions with $p_{R, p}(X) > 0$ in this study, this distribution has smallest value of $p_{R, p}(X)$.

The next two distributions that we consider belong to the class of so-called p-max-stable laws. They can be generalized to a strictly increasing affine transformation and $p_{L,p}$, and $p_{R,p}$ characteristics will not change.
Using power normalizations Pancheva (1985) \cite{Pancheva} obtained them as limiting laws of power transformed maximums. Falk et al. (2004) \cite{Falk2004} describe domains of attraction of these laws under power normalization. We have already seen in Corollary 1 of  e), Theorem 1 that these transformations increase the values of $p_{R, p}$.   Ravi and Saeb (2012) \cite{Ravi} have obtained their entropies.

 \begin{itemize}
\item  $H_{1}$ type. (See Pancheva (1985) \cite{Pancheva}). Let $\alpha > 0$, and $X \in H_1(\alpha)$, i.e.
\end{itemize}

\begin{equation}\label{H1CDF}
  F_X(x) =\left\{\begin{array}{ccc}
                   0 & ,  & x < 1 \\
                   exp\left\{-(log\,\,x)^{-\alpha}\right\} & , & x \geq 1
                 \end{array}
\right..
\end{equation}
 The quantile function of this distribution is
$F_X^\leftarrow(p) = exp\left\{(-log\,\, p)^{-1/\alpha}\right\}$. \footnote{It can be seen e.g. in the supplementary material of Soza et al. (2019) \cite{JordanovaStehlik2018} who consider the case $p = 0.25$.}
Therefore
$$p_{L,p}(X) = P\left\{X < \frac{1}{p}exp\left\{[-ln\,\,(p)]^{-1/\alpha}\right\} - \frac{1-p}{p}exp\left\{[-ln\,\,(1- p)]^{-1/\alpha}\right\}\right\}.$$
and for $exp\left\{[-ln\,\,(p)]^{-1/\alpha}\right\} - (1-p)exp\left\{[-ln\,\,(1- p)]^{-1/\alpha}\right\} > p$, $p_{L,p}(X)$ is equal to
 $$ exp\left\{-\left\{log\left[\frac{1}{p}exp\left\{[-ln\,\,(p)]^{-1/\alpha}\right\} - \frac{1-p}{p}exp\left\{[-ln\,\,(1- p)]^{-1/\alpha}\right\}\right]\right\}^{-\alpha}\right\},$$
 otherwise $p_{L,p}(X) = 0$.  The last mean that for large $\alpha > 0$ it is possible to observe also left outside values.

  Now let us consider the right tail. As far as for all $\alpha > 0$
 $$\frac{1}{p}exp\left\{[-log\,\, (1-p)]^{-1/\alpha}\right\} - \frac{1-p}{p}exp\left\{[-log(p)]^{-1/\alpha}\right\} > 1$$
 $$exp\left\{[-log\,\, (1-p)]^{-1/\alpha}\right\} - (1 - p)exp\left\{[-log(p)]^{-1/\alpha}\right\} > p$$
 $$exp\left\{[-log\,\, (1-p)]^{-1/\alpha}\right\} - exp\left\{[-log(p)]^{-1/\alpha}\right\} > 0$$
 and $$0 > p\left\{1 - exp\left\{[-log(p)]^{-1/\alpha}\right\}\right\}$$
\begin{eqnarray}
\nonumber  p_{R, p}(X) &=& P\left\{X > \frac{1}{p}exp\left\{[-log\,\, (1-p)]^{-1/\alpha}\right\} - \frac{1-p}{p}exp\left\{[-log(p)]^{-1/\alpha}\right\}\right\} =\\
\nonumber  &=& 1 - exp\left\{-\left\{log\left\{\frac{1}{p}exp\left\{[-log\,\, (1-p)]^{-1/\alpha}\right\} \right.\right.\right.\\
   &-& \left. \left. \left. \frac{1-p}{p}exp\left\{[-log(p)]^{-1/\alpha}\right\}\right\} \right\} ^{-\alpha}\right\} \label{peRH1}
\end{eqnarray}
Figure \ref{fig:peRalphaExtremes} shows their dependence  on $\alpha$ in case $p = 0.25$. This case is considered in the supplementary material of Soza et al. (2019) \cite{JordanovaStehlik2018}.

\begin{itemize}
\item  $H_{4}$ type.   Let $\alpha > 0$,
$$F_X(x) =\left\{\begin{array}{ccc}
                   exp\left\{-[ln(-\,x)]^{\alpha}\right\} & , & x < -1\\
                   1 & ,  & x \geq -1
                 \end{array}
\right..$$
\end{itemize}
It is one of the limiting distributions of power-transformed maxima obtained in Pancheva (1984) \cite{Pancheva}. Ravi and Saeb (2012) \cite{Ravi} calculate its Shannon entropy. $H_{4}$ type is known also as log-Weibull law, and it is one of the p-max stable laws.
The quantile function has the form $F_{Y}^\leftarrow(p) = -exp\{[-log(p)]^{1/\alpha}\}$, $p \in (0, 1).$
Therefore
$$p_{L,p}(X) = P\left\{X < \frac{1-p}{p}exp\{[-log(1 - p)]^{1/\alpha}\}-\frac{1}{p}exp\{[-log(p)]^{1/\alpha}\}  \right\}.$$
And for $exp\{[-log(p)]^{1/\alpha}\} - (1-p)exp\{[-log(1 - p)]^{1/\alpha}\} > p$ we have that $p_{L,p}(X)$ is equal to
 $$exp\left\{-\left\{log\left\{\frac{1}{p}exp\{[-log(p)]^{1/\alpha}\} - \frac{1-p}{p}exp\{[-log(1 - p)]^{1/\alpha}\}\right\}\right\}^{\alpha}\right\},$$
and  $p_{L,p}(X) = 0$ otherwise.

  When consider the right tail for all $\alpha > 0$ as far as $p < 1 - p$
 $$\frac{1-p}{p}exp\{[-log(p)]^{1/\alpha}\} - \frac{1}{p}exp\{[-log(1 - p)]^{1/\alpha}\} < - 1$$
 $$(1-p)exp\{[-log(p)]^{1/\alpha}\}-exp\{[-log(1 - p)]^{1/\alpha}\} < -p$$
 $$-(1-p)exp\{[-log(p)]^{1/\alpha}\}-\left\{-exp\{[-log(1 - p)]^{1/\alpha}\}\right\} > $$ $$> -(1-p)exp\{[-log(p)]^{1/\alpha}\}-\left\{-exp\{[-log(p)]^{1/\alpha}\}\right\}$$
 The last expression is equal to $-(-p)exp\{[-log(p)]^{1/\alpha}\} > p$,   therefore $p_{R, p}(X)$ is equal to
\begin{eqnarray*}
    & &  P\left\{X > \frac{1-p}{p}exp\{[-log(p)]^{1/\alpha}\} - \frac{1}{p}exp\{[-log(1 - p)]^{1/\alpha}\}\right\}\\
  &=& 1 - exp\left\{-\left\{log\left\{  \frac{1}{p}exp\{[-log(1 - p)]^{1/\alpha}\} - \frac{1-p}{p}exp\{[-log(p)]^{1/\alpha}\}\right\}\right\}^{\alpha}\right\}
\end{eqnarray*}
As in the previous case for large values of $\alpha > 0$ it is possible to observe both left and right outside values.
The dependance of $p_{R, 0.25}(-X)$ on $\alpha$ is depicted on Figure \ref{fig:peRalphaExtremes}. We call the distribution of $-X$, "$H_{4}$ positive", and we have denoted it by $H_{4}^+$.

Note that the function
 $g(\alpha) = exp\{[log(4)]^{1/\alpha}\} - \frac{3}{4}exp\{[log\frac{4}{3}]^{1/\alpha}\} - \frac{1}{4}$ is decreasing in $\alpha$ and $g(2) \approx 1.713627$ therefore for $\alpha \in (0, 2]$
 \begin{eqnarray*}
   p_{L,0.25}(X) &=& p_{R,0.25}(-X) \\
    &=& exp\left\{-\left\{log\left\{4exp\{[log(4)]^{1/\alpha}\} - 3exp\{[log(\frac{4}{3})]^{1/\alpha}\}\right\}\right\}^{\alpha}\right\}.
 \end{eqnarray*}
See the supplementary material of Soza et al. \cite{JordanovaStehlik2018}.
 \begin{itemize}
\item  $log-Par_\alpha$ type.   Log-Pareto law with parameter $\alpha$ seems to be introduced in Cormann and Reiss (2009) \cite{CormannAndReiss2009}.
 More precisely here we assume that
\end{itemize}
$$F_X(x) =\left\{\begin{array}{ccc}
                   0 & ,  & x < e \\
                   1 -(log(x))^{-\alpha} & , & x \geq e
                 \end{array}
\right..$$

 This distribution belongs to $\Pi$ class considered e.g. in de Haan and Ferreira (2006) \cite{deHaanFerreira} or Embrehts et al. (1997) \cite{EMK}. Ravi and Saeb (2012) \cite{Ravi} investigate their entropies.
 Due to Corollary 1 of Theorem 1 its tail is heavier than the tail of Fr$\acute{e}$chet distribution. The quantile function of this distribution is
$F_X^\leftarrow(p) = exp\left\{(1 -  p)^{-1/\alpha}\right\}$. Therefore
 for $\frac{1}{p}exp\left\{(1 -  p)^{-1/\alpha}\right\} - \frac{1-p}{p}exp\left\{p^{-1/\alpha}\right\} > e$
 $$p_{L,p}(X) = 1 -\left\{log\left[\frac{1}{p}exp\left[(1 -  p)^{-1/\alpha}\right] - \frac{1-p}{p}exp\left(p^{-1/\alpha}\right)\right]\right\}^{-\alpha},$$
 otherwise $p_{L,p}(X) = 0$.

  When consider the right tail for all $\alpha > 0$ as far as $p < 1 - p$
   $$exp\left(p^{-1/\alpha}\right) + p\, exp\left[(1-p)^{-1/\alpha}\right] \geq  exp\left[(1-p)^{-1/\alpha}\right] + ep$$
  $$exp\left(p^{-1/\alpha}\right) - (1-p)exp\left[(1 -  p)^{-1/\alpha}\right] \geq ep$$
 $$\frac{1}{p}exp\left(p^{-1/\alpha}\right) - \frac{1-p}{p}exp\left[(1 -  p)^{-1/\alpha}\right] \geq e$$
   therefore from the definition of $p_{R, p}(X)$ we obtain
\begin{eqnarray*}
  p_{R, p}(X) &=& P\left\{X > \frac{1}{p}exp\left(p^{-1/\alpha}\right) - \frac{1-p}{p}exp\left[(1 -  p)^{-1/\alpha}\right]\right\}\\
  &=& \left\{log\left\{\frac{1}{p}exp\left(p^{-1/\alpha}\right) - \frac{1-p}{p}exp\left[(1 -  p)^{-1/\alpha}\right]\right\}\right\}^{-\alpha}.
\end{eqnarray*}

 The dependence of  $p_{R, 0.25}(X)$ on $\alpha$ could be seen on Figure \ref{fig:peRalphaExtremes}. We observe that its tail behaviour almost coincide with log-Fr$\acute{e}$chet, i.e. $H_1$, and within the considered distributions, for fixed $\alpha$ according to $p_{R, 0.25}(X)$ the last two distributions have highest probabilities to observe extreme outside values.  Having in mind that without transformations the tails of Fr$\acute{e}$chet and Pareto distributions are heavy-tailed Cl. Neves et al. (2008) \cite{ClNeves}, Corman and Reiss (2009) \cite{CormannAndReiss2009} or Falk (2004) \cite{Falk2004} call them "super heavy-tailed".

Further on in this section, we consider distributions which quantile function have no explicit form. Therefore we use R software (2018) \cite{R} in order to obtain obtain $F^\leftarrow(0.25)$ and $F^\leftarrow(0.75)$. Then we come back to the well-known formulas for c.d.f. and obtain $p_{R, 0.25}(X)$ characteristics.

\begin{itemize}
\item  Normal distribution. Assume $\mu \in R$, $\sigma^2 > 0$ and $X \in N(\mu, \sigma^2)$.
    \end{itemize}
   W.l.g. $\mu = 0$ and $\sigma^2 = 1$. Due to the symmetry of this distribution with respect to (w.r.t.) Oy, for all $p \in (0,0.5)$ we have that $p_{L, p}(X) = p_{R, p}(X)$. In particular $p_{L, 0.25}(X) = p_{R, 0.25}(X) \approx 0.000001171$.

\begin{itemize}
\item  $t$-distribution. Assume $n \in N$ and $X \in t(n)$.
\end{itemize}
The symmetry of the p.d.fs. of these distributions w.r.t. $Oy$ implies $p_{L, p}(X) = p_{R, p}(X)$, for all $p \in (0,0.5)$. The values of these characteristics for $n = 1, 2, ..., 10$  are presented in Table \ref{tab:t}, and could be seen also in Jordanova and Petkova \cite{MoniPoli2018}, and in the supplementary material of Soza et al. \cite{JordanovaStehlik2018}.

\begin{center}
\begin{table}
\begin{center}
\begin{tabular}{|c|c|c|c|c|c|}
  \hline
  n & 1 & 2 & 3 & 4 & 5  \\
  \hline
 $p_{L, 0.25}(X) = p_{R, 0.25}(X)$ & 0.0452 & 0.0146 & 0.0064 & 0.0033 & 0.0019  \\
    \hline
    \hline
    n &  6 & 7 & 8 & 9 & 10 \\
     \hline
  $p_{L, 0.25}(X) = p_{R, 0.25}(X)$& 0.0012  & 0.0008 & 0.0006 & 0.0004 & 0.0003 \\
    \hline
\end{tabular}
\caption{The dependence of $p_{L, 0.25}(X) = p_{R, 0.25}(X)$ on $n$. $X \in t(n)$.}\label{tab:t}
\end{center}
\end{table}
\end{center}

\begin{itemize}
\item Gamma distribution. Assume $\alpha > 0$, $\beta > 0$ and $X \in Gamma(\alpha, \beta)$ which means that
\end{itemize}
$$F_X(x) =\left\{\begin{array}{ccc}
                   0 & ,  & x < 0 \\
                   \int_0^x \frac{\beta^\alpha}{\Gamma(\alpha)} y^{\alpha - 1} e^{-\beta y} dy & , & x \geq 0
                 \end{array}
\right..$$
W.l.g. we assume that $\beta = 1$.  The plot of the dependence of  $p_{R, 0.25}(X)$ on $\alpha$ is depicted on Figure \ref{fig:peRalphaE}.

 \begin{itemize}
\item Hill-horror distribution. For $\alpha > 0$  Embrechts et al. (1997) \cite{EMK} define it via its quantile function
\end{itemize}
\begin{equation}\label{HHQuantilefunction}
  F_X^\leftarrow(p) = \frac{-log(1 - p)}{\sqrt[\alpha]{1 - p}}, \quad p \in (0, 1).
\end{equation}

Then $p_{L, 0.25}(X) = 0$. For the values of $p_{R, 0.25}(X)$ see Figure \ref{fig:peRalphaE}.  We observe that this distribution has one of the heaviest right tails within the considered probability types.






Following this approach, we can find explicit values or plots of $p_{L, p}$ and $p_{R, p}$ characteristics for many other distributions and in this way to compare their tails. For example log-Positive Weibull, log-Gumbel, Invence-Gamma, Log-Gamma, Beta prime,  or powers bigger than one, of these and other distributions.

\section{Properties of the estimators}
\label{sec:3}
In the previous section, we have considered some particular cases of probability laws and we have shown which parameter governs the heaviness of the tail of the corresponding distribution according to our classification. For the distributions with regularly varying tails, it coincides with the very well-known index of regular variation. In this section, we obtain different asymptotic properties of the estimators of the corresponding parameters of heaviness of the tails. The general formula for the joint distribution of order statistics is very well known. Together with the formula for their conditional distributions they could be found e.g. in Nevzorov (2001) \cite{Nevzorov} or in Arnold et al. (1992) \cite{Arnold1992}. The following lemma is their immediate corollary. Its first part summarises the same results in the terms of equality in distributions. In vi) we have expressed the bivariate vector of order statistics as a bivariate function of independent r.vs. This allows as to make the same with the fences in the next property. Finally two explicit formulae for the probability mass functions of the numbers of left and right $p$-outside values in a sample of independent observations are presented.  Due to their complicated forms further on the section proceeds with asymptotic results.

\begin{lm}\label{lem:lm1} If $\xi \in Beta(i, j)$, $i \in N$, $j \in N$, $F(z) = P(X < z)$ is some c.d.f. of a r.v. $X$ and $c > 0$ is a constant, then

\begin{enumerate}
\item[i)] $-\frac{log\, (\xi)}{c} \stackrel{d}{=} \epsilon_{(j, i+j-1)},$
where $\epsilon_{(j, i+j-1)}$ is the $j$-th order statistics in a sample of $i+j-1$ independent observations on i.i.d. Exponential r.vs. with parameter $c > 0$.
\item[ii)] $\frac{1}{c\xi^\alpha} \stackrel{d}{=} \nu_{(j, i+j-1)}$,
where $\nu_{(j, i+j-1)}$ is the $j$-th order statistics in a sample of $i+j-1$ independent observations on i.i.d. Pareto distributed r.vs. with parameters $\frac{1}{\alpha} > 0$ and $\delta = \frac{1}{c}$.
\item[iii)] $\frac{-1}{c \,\, log(\xi)} \stackrel{d}{=} \phi_{(j, i+j-1)}$, where $\phi_{(j, i+j-1)}$ is the $j$-th order statistics in a sample of $i+j-1$ independent observations on i.i.d. Fr$\acute{e}$chet distributed r.vs. with parameter $\alpha = 1$ and scale parameter $c$.
\item[iv)] $F^\leftarrow(\xi)\stackrel{d}{=} \kappa_{(i, i+j-1)},$ and $(1 - F)^\leftarrow(\xi)\stackrel{d}{=} \kappa_{(j, i+j-1)}$, where $\kappa_{(s, i+j-1)}$ is the $s$-th order statistic of a sample of $i+j-1$ independent observations on a r.v. with absolutely continuous c.d.f. $F$, $s = i, j$.
\item[v)] For $1 \leq i < j$
$$(X_{(j, i+j-1)}|X_{(i, i+j-1)} = x) \stackrel{d}{=} \tilde{\tilde{\theta}}_{(j-i, j-1)} \stackrel{d}{=} G_x^\leftarrow(\xi^*),$$
 $$(X_{(i, i+j-1)}|X_{(j, i+j-1)} = x) \stackrel{d}{=} \tilde{\theta}_{(i, j-1)} \stackrel{d}{=} T_x^\leftarrow(\xi^{**}),$$ where $\tilde{\tilde{\theta}}_{(j-i, j-1)}$ is the $j-i$-th order statistics in a sample of $j-1$ independent observations on i.i.d. r.vs. with $X_{LT}$ c.d.f.
$G_x(y) = \frac{F(y) - F(x)}{1-F(x)}$, $y > x$, $\tilde{\theta}_{(i, j-1)}$ is the $i$-th order statistics in a sample of $j-1$ independent observations on i.i.d. r.vs. with $X_{RT}$ c.d.f.
$T_x(y) = \frac{F(y)}{F(x)}$, $y < x$, $\xi^* \in Beta(j - i, i)$, and $\xi^{**} \in Beta(i, j - i)$.
\item[vi)]   Assume $\xi$ and $\xi^*$ are independent, and $\xi^* \in Beta(j-i, i)$. Denote by $\xi^{**} = 1 - \xi^*$, and $\xi^{***} = 1 - \xi$.\footnote{Then $1 - \xi^* =: \xi^{**} \in Beta(i, j-i)$ and $1 - \xi =: \xi^{***} \in Beta(j, i)$.} Then for $n = i + j - 1$, $i, j \in N$, $1 \leq i < j$,
\begin{enumerate}
\item $(X_{(j, i+j-1)}, X_{(i, i+j-1)})$
\begin{eqnarray*}
   &\stackrel{d}{=}& \{ F^\leftarrow(1 - \xi^{**} \xi^{***}),F^\leftarrow (1 - \xi^{***})\}\\
   &\stackrel{d}{=}& \left\{(1-F)^\leftarrow(\xi^{**} \xi^{***}), (1-F)^\leftarrow(\xi^{***})\right\}\\
   &\stackrel{d}{=}& \left\{\left(\frac{1}{1-F}\right)^\leftarrow\left(\frac{1}{\xi^{**} \xi^{***}}\right), \left(\frac{1}{1-F}\right)^\leftarrow\left(\frac{1}{\xi^{***}}\right)\right\}
\end{eqnarray*}
Moreover $\xi^{**} \xi^{***} \stackrel{d}{=} \xi$.
\item The empirical $\frac{i}{i+j}$-right fences
\begin{eqnarray} \label{rightfences}
 R_n(X,\frac{i}{i+j}) &=& \left(1+\frac{j}{i}\right)X_{j, i + j - 1} -  \frac{j}{i} X_{i, i + j - 1} \\
  \nonumber  &\stackrel{d}{=}& \left(1+\frac{j}{i}\right)(1 - F)^\leftarrow(\xi^{**}\xi^{***}) -  \frac{j}{i}(1- F)^\leftarrow (\xi^{***}).
\end{eqnarray}
 \item The empirical $\frac{i}{i+j}$-left fences
\begin{eqnarray} \label{leftfences}
L_n(X,\frac{i}{i+j}) &=&  \left(1+\frac{j}{i}\right)X_{i, i + j - 1} -  \frac{j}{i} X_{j, i + j - 1} \\
\nonumber &\stackrel{d}{=}& \left(1+\frac{j}{i}\right)(1- F)^\leftarrow (\xi^{***}) -  \frac{j}{i}(1 - F)^\leftarrow(\xi^{**}\xi^{***}).
\end{eqnarray}
 \end{enumerate}
 \item[vii)] For $k = 0, 1, ..., i+j$,

$P(n_R(\frac{i}{i+j},i+j-1) = k)$
\begin{eqnarray*}
    &=& \frac{(i+j-1)!}{k!(i+j-k-1)!}\\
   &.&\int_0^1\int_0^1\left\{\frac{1-F\left[\left(1+\frac{j}{i}\right)(1 - F)^\leftarrow(xy) -  \frac{j}{i}(1 - F)^\leftarrow (x)\right]}{F\left[\left(1+\frac{j}{i}\right)(1 - F)^\leftarrow(xy) -  \frac{j}{i}(1 - F)^\leftarrow (x)\right]}\right\}^k \\
   &.&\frac{x^{j-1}(1-x)^{i-1}}{B(i,j)}\frac{y^{i-1}(1-y)^{j-i-1}}{B(j-i,i)}\\
   &.&\left\{F\left[\left(1+\frac{j}{i}\right)(1 - F)^\leftarrow(xy) -  \frac{j}{i}(1 - F)^\leftarrow (x)\right]\right\}^{i+j-1}dydx
\end{eqnarray*}
$P(n_L(\frac{i}{i+j},i+j-1) = k)$
\begin{eqnarray*}
   &=& \frac{(i+j-1)!}{k!(i+j-k-1)!}\\
   &.&\int_0^1\int_0^1\left\{\frac{F\left[\left(1+\frac{j}{i}\right)(1 - F)^\leftarrow(x) -  \frac{j}{i}(1 - F)^\leftarrow(xy)\right]}{1-F\left[\left(1+\frac{j}{i}\right)(1 - F)^\leftarrow(x) -  \frac{j}{i}(1 - F)^\leftarrow (xy)\right]}\right\}^k \\
   &.&\frac{x^{j-1}(1-x)^{i-1}}{B(i,j)}\frac{y^{i-1}(1-y)^{j-i-1}}{B(j-i,i)}\\
   &.&\left\{1-F\left[\left(1+\frac{j}{i}\right)(1 - F)^\leftarrow(x) -  \frac{j}{i}(1 - F)^\leftarrow (xy)\right]\right\}^{i+j-1}dydx
\end{eqnarray*}

    \end{enumerate}
    \end{lm}

\begin{rem} As an additional result, we can use the above theorem to obtain different univariate and bivariate distributions and new relations between them. This approach is analogous to the one applied by Eugene et al. (2002) \cite{Eugene} or Cordeiro et al. (2012) \cite{Cordeiro} among others, who consider Generalized-Beta generalized distributions.
\end{rem}

\begin{rem} Using the general formula for the moments of order statistics, for $m = 1, 2, ..., n$, and $r \in R$
$$E[X_{(m,n)}^r] = \frac{n!}{(m-1)!(n-m)!}\int_{-\infty}^\infty x^r[F(x)]^{m-1}[1-F(x)]^{n-m}dF(x),$$
which could be seen e.g. in the books of Arnold et al. (1992)\cite{Arnold1992} or Nevzorov (2001) \cite{Nevzorov}, we can easily obtain the general formulae for the mean and the variance of $L_n(X,\frac{i}{i+j})$ and $R_n(X,\frac{i}{i+j})$ in cases when they exist. For example in Pareto (\ref{Pareto}) case, for $i = k$, $j = sk$, $s = 2, 3, ...,$ $L_{k(s+1)-1}(X,\frac{1}{1+s})$ and $R_{k(s+1)-1}(X,\frac{1}{1+s})$, are asymptotically unbiased estimators correspondingly for $L(F, \frac{1}{s + 1})$ and $R(F, \frac{1}{s + 1})$. More precisely
$$\lim_{k\to \infty} EL_{k(s+1)-1}(X,\frac{1}{s+1}) = L(F, \frac{1}{s + 1}), \quad and$$
$$\lim_{k\to \infty} ER_{k(s+1)-1}(X,\frac{1}{s+1}) = R(F, \frac{1}{s + 1}).$$
\end{rem}

The following result is an immediate corollary of the definition of convergence in probability, quantile transform, a.s. convergence of empirical quantiles to the corresponding theoretical one, and Slutsky's theorem about continuous functions. See e.g. Embrechts et al. (1997) \cite{EMK}.

\begin{thm}\label{thm:Lemma2} Given a sample of independent observations, for any fixed $s = 2, 3, ...$
\begin{eqnarray}
  L_{(s+1)k - 1}\left(F, \frac{1}{s + 1}\right) &\stackrel{P}{\to }& L(F, \frac{1}{s + 1}), \quad k \to \infty,\\ \label{LFencesConsistency}
  R_{(s+1)k - 1}\left(F, \frac{1}{s + 1}\right) &\stackrel{P}{\to }& R(F, \frac{1}{s + 1}),  \quad k \to \infty. \label{RFencesConsistency}
\end{eqnarray}
\end{thm}

Cadwell (1953) \cite{Cadwell} finds the distribution of quasi-ranges in samples from a normal population. He gives us the idea about the next result. Rider  (1959) \cite{Rider} obtains their exact distribution in case of samples from an exponential population. Sarhan et al. (1963) \cite{SarhanExponential} propose  simplified estimates in this case. The asymptotic normality of the appropriately normalized univariate distributions of the central order statistics is investigated in Smirnov (1949) \cite{Smirnov1949}. Note that in the next theorem, because of the special choice of the numbers of order statistics, $p$ and $n$ his conditions $\frac{k}{n} \to p \in(0, 1)$  and $\sqrt{n}(\frac{k}{n} - p) \to \mu \in (-\infty, \infty)$ are satisfied. Moreover, in our case $\mu = 0$. The theorem about the joint distribution of the central order statistics of i.i.d. observations, could be seen e.g. in Nair (2013) \cite{Nair2013}, p.330, or Arnold et al. (1992) \cite{Arnold1992}, p. 226, among others. The multivariate delta method is a very powerful technique for obtaining confidence intervals in such cases. It can be seen e.g. in Sobel (1982) \cite{MultivariateDeltaMethod}. In the next theorem we use these results and obtain the limiting distribution of the fences of central order statistics.

\begin{thm}\label{thm:ThmGeneralAsumptoticNormalityFences} Consider a sample of $n = (s+1)k-1$, $s = 1, 2, 3, ...$ observations on a r.v. $X$ with c.d.f. $F$ and p.d.f. $f = F'$. Suppose that there exists $c_{F,s} : = f\left[F^\leftarrow(\frac{1}{s+1})\right] \in (0,  \infty)$ and $d_{F,s} =  f\left[F^\leftarrow(\frac{s}{s+1})\right] \in (0, \infty)$. Then
$$\lim_{k \to \infty} E L_{n}\left(F, \frac{1}{s + 1}\right) = L(F, \frac{1}{s + 1})$$
$$\lim_{k \to \infty} E R_{n}\left(F, \frac{1}{s + 1}\right) = R(F, \frac{1}{s + 1})$$
and for $k \to \infty$
 \begin{equation}\label{General_CI_L_fence}
\sqrt{(s+1)k-1}\left|L_{n}\left(F, \frac{1}{s + 1}\right) - L(F, \frac{1}{s + 1})\right| \stackrel{d}{\to } N\left(0; V_{L,F,s}\right),
 \end{equation}
 \begin{equation}\label{General_CI_R_fence}
\sqrt{(s+1)k-1}\left|R_{n}\left(F, \frac{1}{s + 1}\right) - R(F, \frac{1}{s + 1})\right| \stackrel{d}{\to} N\left(0; V_{R,F,s}\right),
 \end{equation}
  where $V_{L,F,s} = \frac{s}{(s+1)^2} \left[\frac{(s + 1)^2}{c_{F,s}^2} - \frac{2(s+1)}{c_{F,s}d_{F,s}} + \frac{s^2}{d_{F,s}^2}\right]$ and \\$V_{R,F,s} = \frac{s}{(s+1)^2} \left[\frac{s^2}{c_{F,s}^2} - \frac{2(s+1)}{c_{F,s}d_{F,s}} + \frac{(s + 1)^2}{d_{F,s}^2}\right].$
\end{thm}



These allows as to compute the asymptotic confidence intervals of these estimators for $k \to \infty$, $n = (s+1)k - 1$, and  $s = 1, 2, 3, ...$, fixed. Denote $L_{F,s,k} = L_{(s+1)k - 1}\left(F, \frac{1}{s + 1}\right)$ and $R_{F,s,k} = R_{(s+1)k - 1}\left(F, \frac{1}{s + 1}\right)$. If the conditions of Theorem \ref{thm:ThmGeneralAsumptoticNormalityFences} are satisfied, then given $\alpha \in (0, 1)$,
\begin{equation}\label{CI_L}
L_{F,s,k} - z_{\alpha} \sqrt{\frac{V_{L,F,s}}{(s+1)k-1}} \leq L(F, \frac{1}{s + 1}) \leq L_{F,s,k} + z_{\alpha} \sqrt{\frac{V_{L,F,s}}{(s+1)k-1}},
\end{equation}
\begin{equation}\label{CI_R}
R_{F,s,k} - z_{\alpha} \sqrt{\frac{V_{R,F,s}}{(s+1)k-1}} \leq R(F, \frac{1}{s + 1}) \leq R_{F,s,k} + z_{\alpha} \sqrt{\frac{V_{R,F,s}}{(s+1)k-1}}.
\end{equation}

The next theorem explains why different probabilities for outside values can be useful for estimating the tail behaviour of the observed distribution.
For a fixed $s = 2, 3, ...$ and $k \to \infty$, we apply the approach of Dembinska (2012) \cite{Dembinska2012NZ}, for the bivariate case, and obtain that
$\frac{n_L\left(\frac{1}{s+1}; (s+1)k-1\right)}{(s + 1)k - 1}$ and $\frac{n_R\left(\frac{1}{s+1};(s+1)k-1\right)}{ (s + 1)k - 1}$
are strongly consistent estimators correspondingly of $p_{L,\frac{1}{s+1}}(X)$ and $p_{R,\frac{1}{s+1}}(X)$. Moreover Dembinska (2017) \cite{Dembinska2017Metrika} shows that this approach works not only for i.i.d., but also also for  strictly stationary and ergodic sequences.
\begin{thm}\label{thm:Strongconsistencofpel} Let $s = 2, 3, ...$ be fixed. Assume $f[F^\leftarrow(\frac{1}{s+1})] \in (0, \infty)$ and \\ $f[F^\leftarrow(\frac{s}{s+1})] \in (0, \infty)$.
\begin{enumerate}
  \item If $f[F^\leftarrow(\frac{1}{s+1})] \in (0, \infty)$, then  $$\frac{n_L\left(\frac{1}{s+1}; (s+1)k-1\right)}{(s + 1)k - 1} {\mathop{\to}\limits_{k \to \infty}^{a.s.}} p_{L,\frac{1}{s+1}}(X).$$
  \item If $f[F^\leftarrow(\frac{1}{s+1})] \in (0, \infty)$, then  $$\frac{n_R\left(\frac{1}{s+1};(s+1)k-1\right)}{ (s + 1)k - 1} {\mathop{\to}\limits_{k \to \infty}^{a.s.}} p_{R,\frac{1}{s+1}}(X).$$
\end{enumerate}
\end{thm}

\section{Simulation study}
\label{sec:sym}

In this section we assume that $\mathbf{X}_1, \mathbf{X}_2, ..., \mathbf{X}_n$ are independent realizations of $\mathbf{X}$, with c.d.f. $F$.
Jordanova and Petkova (2017-2018) \cite{MoniPoli2017,MoniPoli2018} assume that $F$ has regularly varying tail. More precisely they consider only the cases when there exists $\alpha > 0$, such that for all $x > 0$,
\begin{equation}\label{RV}
  \lim_{t \to \infty} \frac{1 - F(xt)}{1 - F(t)} = x^{-\alpha}.
\end{equation}
The number $-\alpha$ is called "index of regular variation of the tail of c.d.f.", see e.g. de Haan and Ferreira (2006) \cite{deHaanFerreira} or Resnick (1987) \cite{Resnick87}. Using the explicit form of the corresponding probabilities for extreme outliers according the definition in Devore(2015) \cite{Devore},($0.25$-outliers) Jordanova and Petkova (2017-2018) \cite{MoniPoli2017,MoniPoli2018} obtain distribution sensitive estimators of the unknown parameter $\alpha$ which governs the tail of the considered distributional type. The algorithm consists of the following three main steps.
\begin{enumerate}
  \item Using the results from the previous two sections the explorer chooses the most appropriate probability type (let us call it $T$) for modeling the tail of the distribution of the observed r.v.
  \item Using the formula for $p_{R, 0.25}$ in case $T$ one expresses the unknown parameter $\alpha$.
  \item Replace the theoretical characteristics in the previous step with the corresponding estimators and obtain a new estimator for the parameter which governs the tail behavior.
\end{enumerate}
In their work Jordanova and Petkova (2017-2018) \cite{MoniPoli2017,MoniPoli2018} compare the obtained in this way estimators in Pareto, Fr$\acute{e}$chet, and Hill-Horror case with Hill, t-Hill, Pickands, and Deckers-Einmahl-de Haan estimators and depict the results via a simulation study.  Here we consider two more cases: Log-Logistic case (\ref{LolLogidticCDF}) and $H_1$ case (\ref{H1CDF}). Although in the last case the right tail of the c.d.f. is not regularly varying the next study shows that the approach still gives very good results.

In any of the following five examples using the functions implemented in R (2018), \cite{R} we have simulated $m = 1000$ samples of $n$ independent observations separately on $X$. Then for any fixed $n = 10, 11, ..., 500$ and for any fixed sample we have computed the estimators $\hat{p}_R(0.25, n)= \frac{n_{R}(0.25, n)}{n},$ $F_n^\leftarrow(0.25)$, $F_n^\leftarrow(0.75)$, and $\hat{\alpha}_{\bullet,n}$. Here $n_{R}(0.25, n)$ is the numbers of right extreme outside values in the considered sample of $n$ independent observations, and $\bullet$ means  one of the abbreviations $Par$, $Fr$, $HH$, $H_1$ or $LL$ explained below. Finally we have fixed one of the last estimators and we have averaged the corresponding values of $\hat{\alpha}_{\bullet,n}$ over the considered $n$. The next Figures \ref{fig:SymStudy1}-\ref{fig:SymStudy10} depict the dependence of these values, together with the corresponding asymptotic normal 95\%  confidence intervals,  on the real type of the observed r.v., and on the sample size, for $\alpha = 0.5$, or $1$. We have chosen only the cases when $\alpha$ is small because our observations show that for a fixed sample size the more the outside values, the heavier of the tail of the c.d.f. is and the better the corresponding estimator is.

Let us depict this approach with some examples. In any of them we suppose that $F_n^\leftarrow\left(\frac{3}{4}\right) > 1$ and $\hat{p}_R(0.25, n) > 0$.

\begin{ex} Assime $\mathbf{X} \in Par(\alpha, \delta)$. See (\ref{Pareto}). Having a sample of $n$ independent observations on $\mathbf{X}$, analogously to the generalized method of moments, and following the above algorithm Jordanova and Petkova (2018) \cite{MoniPoli2018} obtain
\begin{equation}\label{EstParn}
\hat{\alpha}_{Par,n} = -\frac{log\,\,\hat{p}_R(0.25, n)}{log\left\{F_n^\leftarrow\left(\frac{3}{4}\right) + 3\left[F_n^\leftarrow\left(\frac{3}{4}\right) - F_n^\leftarrow\left(\frac{1}{4}\right)\right]\right\}}.
\end{equation}
\end{ex}

\begin{ex} If $\mathbf{X} \in Fr(\alpha, c)$, see (\ref{Frechet}), Jordanova and Petkova (2018) \cite{MoniPoli2018} propose
\begin{equation}\label{EstFr}
\hat{\alpha}_{Fr,n} = -\frac{log\{-log[1-\hat{p}_R(0.25, n)]\}}{log\left\{F_n^\leftarrow\left(\frac{3}{4}\right) + 3\left[F_n^\leftarrow\left(\frac{3}{4}\right)- F_n^\leftarrow\left(\frac{1}{4}\right)\right]\right\}}.
\end{equation}
\end{ex}

\begin{ex} Let $\mathbf{X}$ be Hill-Horror distributed. This distribution is usually defined via its quantile function (\ref{HHQuantilefunction}).  Given $F_n^\leftarrow\left(\frac{3}{4}\right) + 3\left[F_n^\leftarrow\left(\frac{3}{4}\right)- F_n^\leftarrow\left(\frac{3}{4}\right)\right] \not= -log\, \hat{p}_R(0.25, n)$ Jordanova and Petkova (2018) \cite{MoniPoli2018} use
\begin{equation}\label{EstHH}
\hat{\alpha}_{HH,n} = \frac{log\, \hat{p}_R(0.25, n)}{log\,\left\{\frac{-log\,\hat{p}_R(0.25, n)}{F_n^\leftarrow\left(\frac{3}{4}\right) + 3\left[F_n^\leftarrow\left(\frac{3}{4}\right)- F_n^\leftarrow\left(\frac{1}{4}\right)\right]}\right\}}.
\end{equation}
\end{ex}

The next two estimators seems to be new. They show that this approach can be applied in much wider than the regularly varying case.

\begin{ex} Let $\mathbf{X} \in H_1$, see (\ref{H1CDF}). It is difficult to solve (\ref{peRH1}) with respect to $\alpha$, therefore we solve the equation
$$ p_{R, 0.25}(X) = 1 - exp\{-[log\,R(X, 0.25)]^{-\alpha}\}.$$
When  express $\alpha$ and replace the theoretical characteristics with the corresponding empirical one we obtain the estimator
\begin{equation}\label{EstH1}
\hat{\alpha}_{H_1,n} = -\frac{log\{-log[1-\hat{p}_R(0.25, n)]\}}{log\left\{log\left\{F_n^\leftarrow\left(\frac{3}{4}\right) + 3\left[F_n^\leftarrow\left(\frac{3}{4}\right) - F_n^\leftarrow\left(\frac{1}{4}\right)\right]\right\}\right\}}.
\end{equation}
\end{ex}

\begin{ex} Suppose $\mathbf{X}$ follow Log-logistic probability law (\ref{LolLogidticCDF}). The equation (\ref{pRpLogLogistic}) have no explicit solution for $\alpha$, therefore we solve the equation
$$ p_{R, 0.25}(X) = \frac{1}{1 + R(X, 0.25)^{\alpha}}.$$
Then we replace the theoretical characteristics with the corresponding empirical one we obtain the estimator
\begin{equation}\label{EstLL}
\hat{\alpha}_{LL,n} = \frac{log\left[\frac{1}{\hat{p}_R(0.25, n)}-1\right]}{log\left\{F_n^\leftarrow\left(\frac{3}{4}\right) + 3\left[F_n^\leftarrow\left(\frac{3}{4}\right) - F_n^\leftarrow\left(\frac{1}{4}\right)\right]\right\}}.
\end{equation}
\end{ex}

Figures \ref{fig:SymStudy1}-\ref{fig:SymStudy10} depict the dependence of these estimators, together with their empirical 95\%  confidence intervals  on the sample size, probability law of the simulated r.v., and the estimated parameter $\alpha$. The names of the estimators in these figures are abbreviated as follows: $\hat{\alpha}_{Par,n} = aParn$, $\hat{\alpha}_{Fr,n} = aFrn$, $\hat{\alpha}_{HH,n} = aHHn$, $\hat{\alpha}_{LL,n} = aLLn$, and $\hat{\alpha}_{H_1,n} = aH1n$.

\begin{figure}
\begin{minipage}[t]{0.5\linewidth}
\includegraphics[scale=.47]{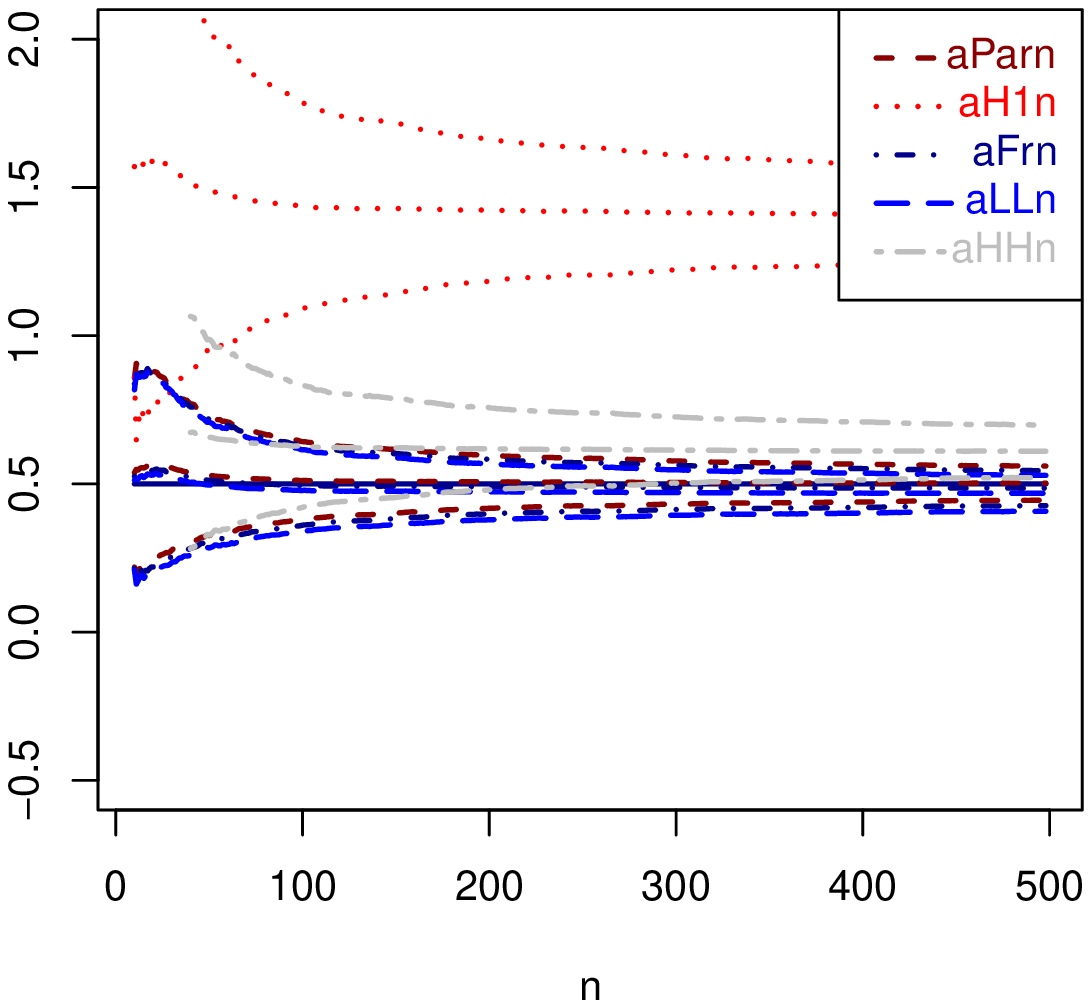}\vspace{-0.3cm}
\caption{Pareto case, $\alpha = 0.5$: Dependence of $\hat{\alpha}_{Par,n} $, $\hat{\alpha}_{Fr,n} $, $\hat{\alpha}_{HH,n}$, $\hat{\alpha}_{H_1,n}$, $\hat{\alpha}_{LL,n}$ on the sample size}
    \label{fig:SymStudy1}
    \end{minipage}
\begin{minipage}[t]{0.5\linewidth}
\includegraphics[scale=.47]{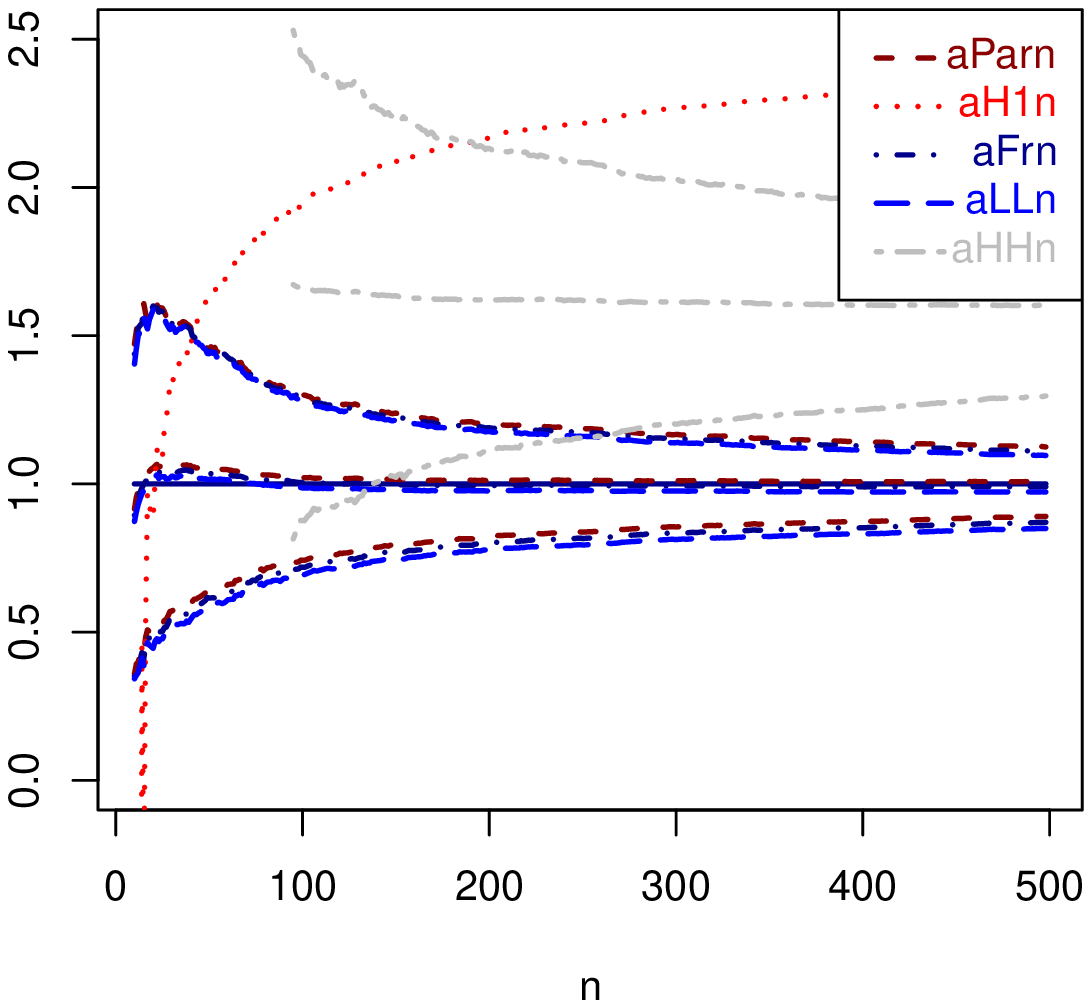}\vspace{-0.3cm}
\caption{Pareto case, $\alpha = 1$: Dependence of $\hat{\alpha}_{Par,n} $, $\hat{\alpha}_{Fr,n} $, $\hat{\alpha}_{HH,n}$, $\hat{\alpha}_{H_1,n}$, $\hat{\alpha}_{LL,n}$ on the sample size}
    \label{fig:SymStudy2}
\end{minipage}
\end{figure}

\begin{figure}
\begin{minipage}[t]{0.5\linewidth}
\includegraphics[scale=.47]{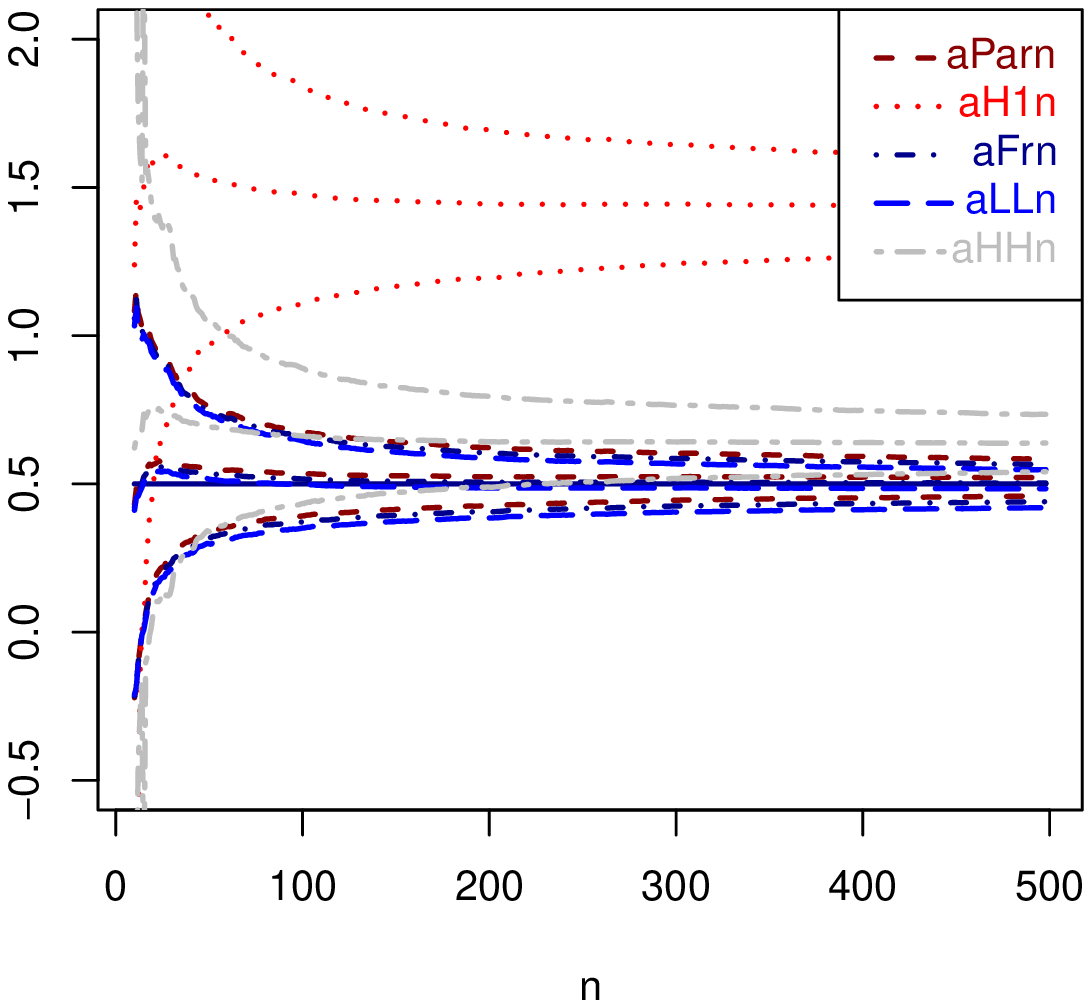}\vspace{-0.3cm}
\caption{Fr$\acute{e}$chet case, $\alpha = 0.5$: Dependence of $\hat{\alpha}_{Par,n} $, $\hat{\alpha}_{Fr,n} $, $\hat{\alpha}_{HH,n}$, $\hat{\alpha}_{H_1,n}$, $\hat{\alpha}_{LL,n}$ on the sample size}
    \label{fig:SymStudy3}
    \end{minipage}
\begin{minipage}[t]{0.5\linewidth}
\includegraphics[scale=.47]{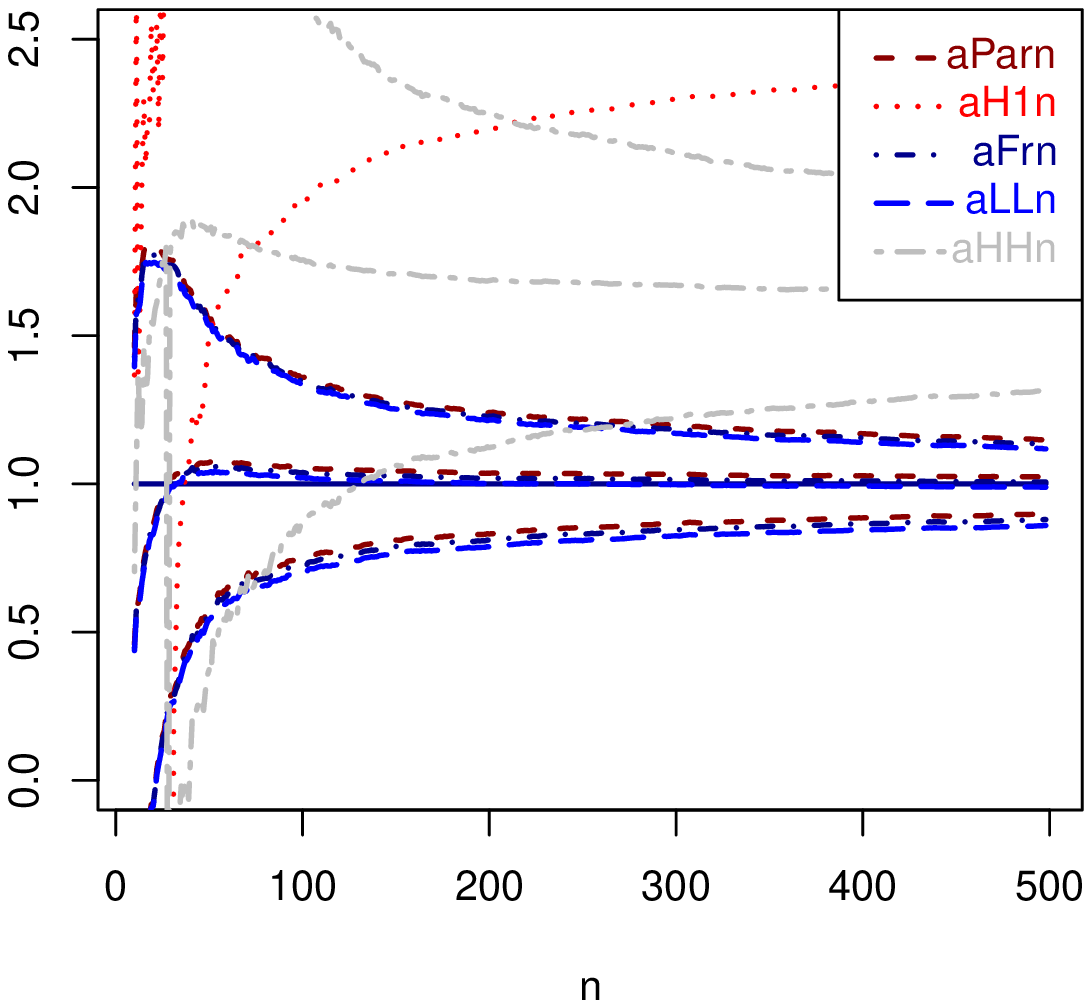}\vspace{-0.3cm}
\caption{Fr$\acute{e}$chet case, $\alpha = 1$: Dependence of $\hat{\alpha}_{Par,n} $, $\hat{\alpha}_{Fr,n} $, $\hat{\alpha}_{HH,n}$, $\hat{\alpha}_{H_1,n}$, $\hat{\alpha}_{LL,n}$ on the sample size}
    \label{fig:SymStudy4}
\end{minipage}
\end{figure}

\begin{figure}
\begin{minipage}[t]{0.5\linewidth}
\includegraphics[scale=.47]{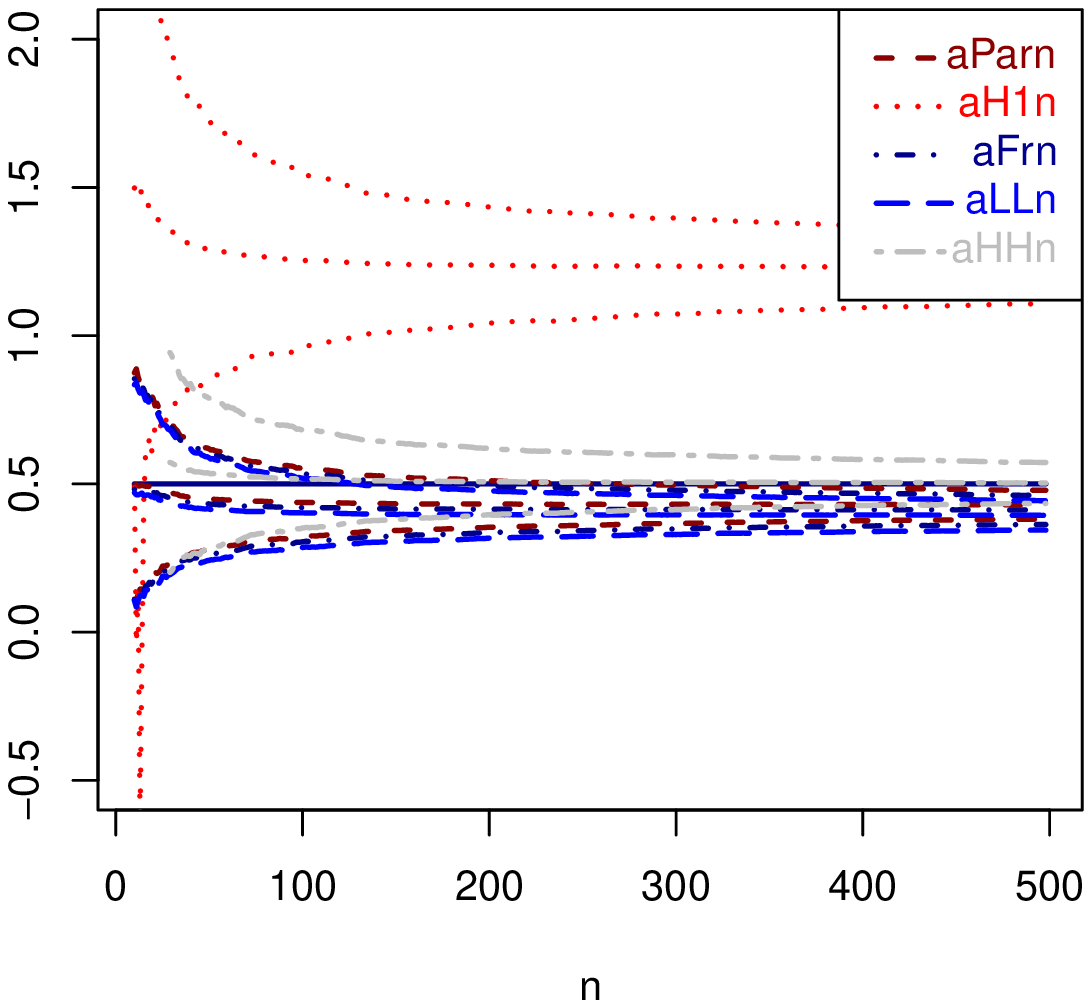}\vspace{-0.3cm}
\caption{Hill-Horror case, $\alpha = 0.5$: Dependence of $\hat{\alpha}_{Par,n} $, $\hat{\alpha}_{Fr,n} $, $\hat{\alpha}_{HH,n}$, $\hat{\alpha}_{H_1,n}$, $\hat{\alpha}_{LL,n}$ on the sample size}
    \label{fig:SymStudy5}
    \end{minipage}
\begin{minipage}[t]{0.5\linewidth}
\includegraphics[scale=.47]{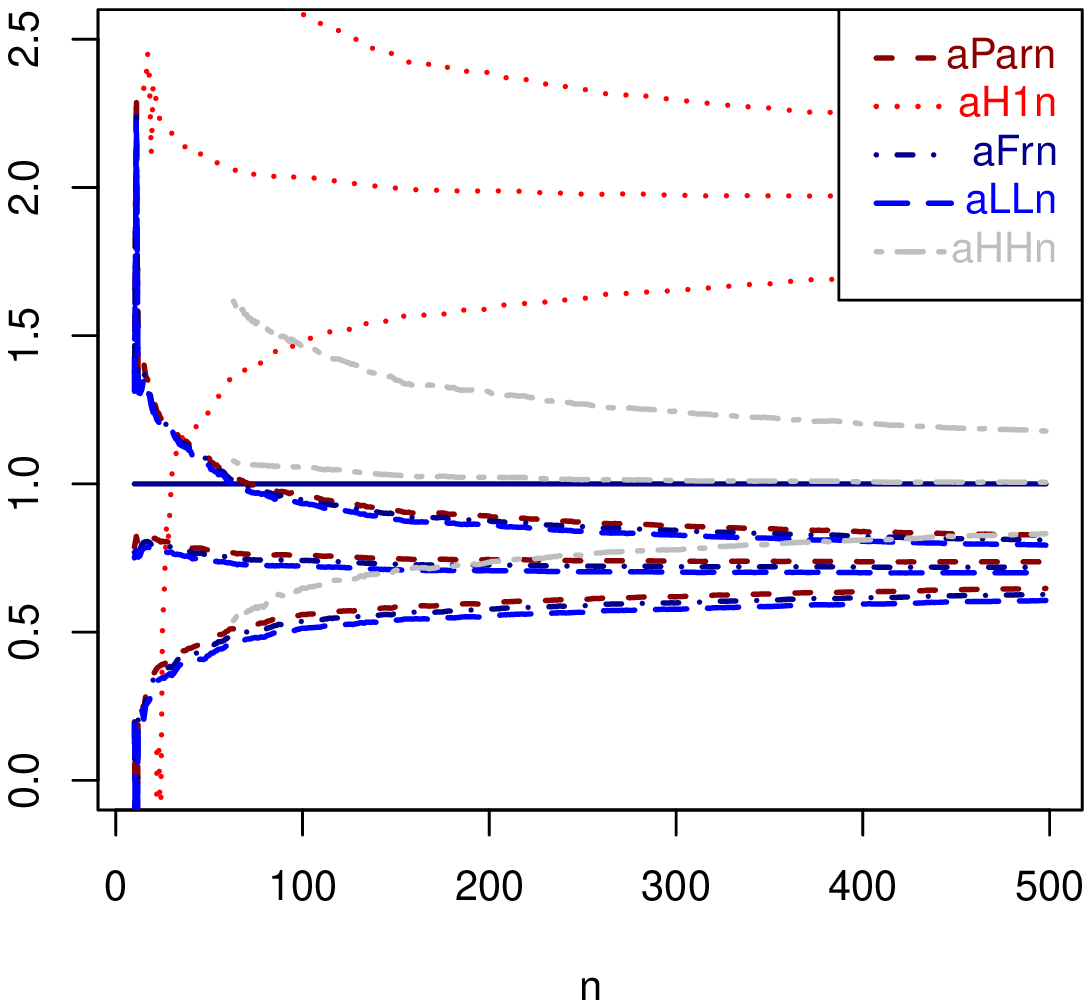}\vspace{-0.3cm}
\caption{Hill-Horror case, $\alpha = 1$: Dependence of $\hat{\alpha}_{Par,n} $, $\hat{\alpha}_{Fr,n} $, $\hat{\alpha}_{HH,n}$, $\hat{\alpha}_{H_1,n}$, $\hat{\alpha}_{LL,n}$ on the sample size}
    \label{fig:SymStudy6}
\end{minipage}
\end{figure}

\begin{figure}
\begin{minipage}[t]{0.5\linewidth}
\includegraphics[scale=.47]{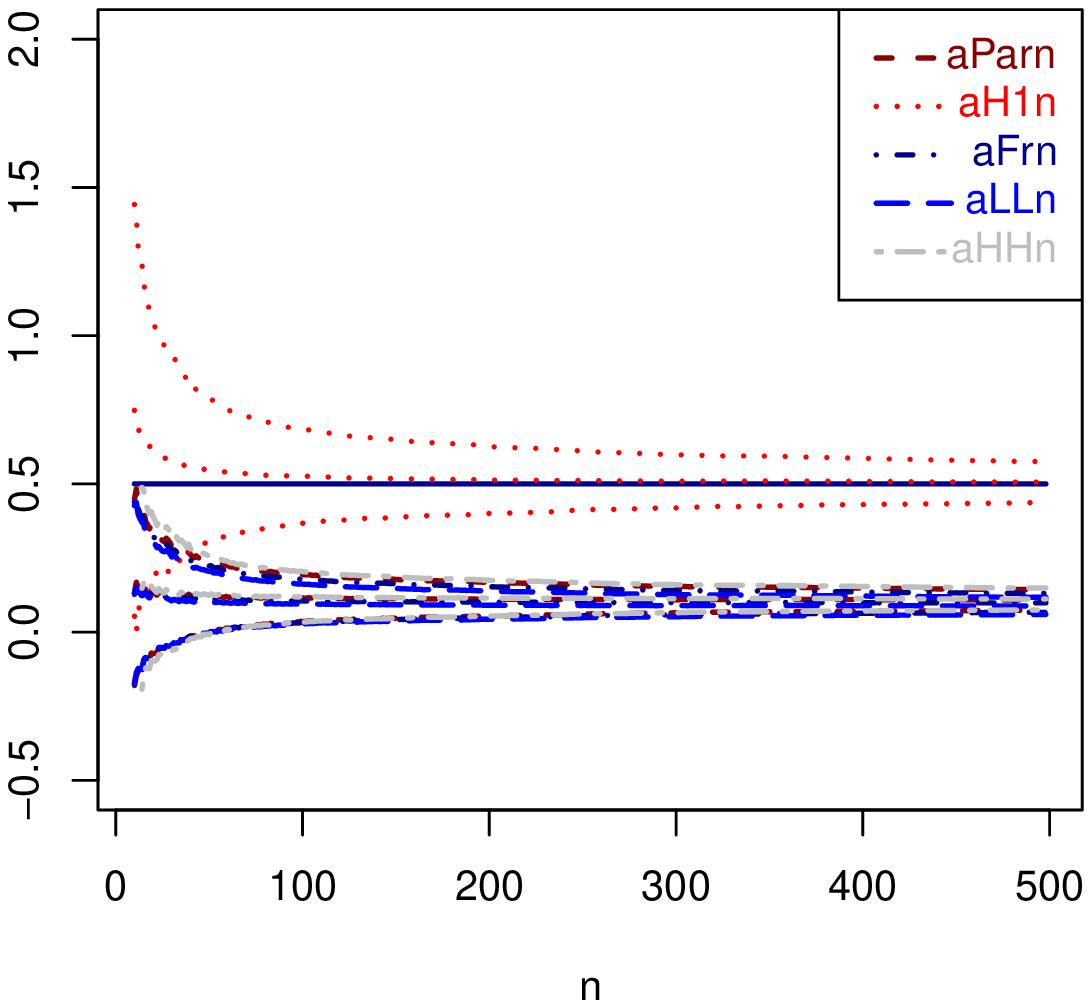}\vspace{-0.3cm}
\caption{$H_1$ case, $\alpha = 0.5$: Dependence of $\hat{\alpha}_{Par,n} $, $\hat{\alpha}_{Fr,n} $, $\hat{\alpha}_{HH,n}$, $\hat{\alpha}_{H_1,n}$, $\hat{\alpha}_{LL,n}$ on the sample size}
    \label{fig:SymStudy7}
    \end{minipage}
\begin{minipage}[t]{0.5\linewidth}
\includegraphics[scale=.47]{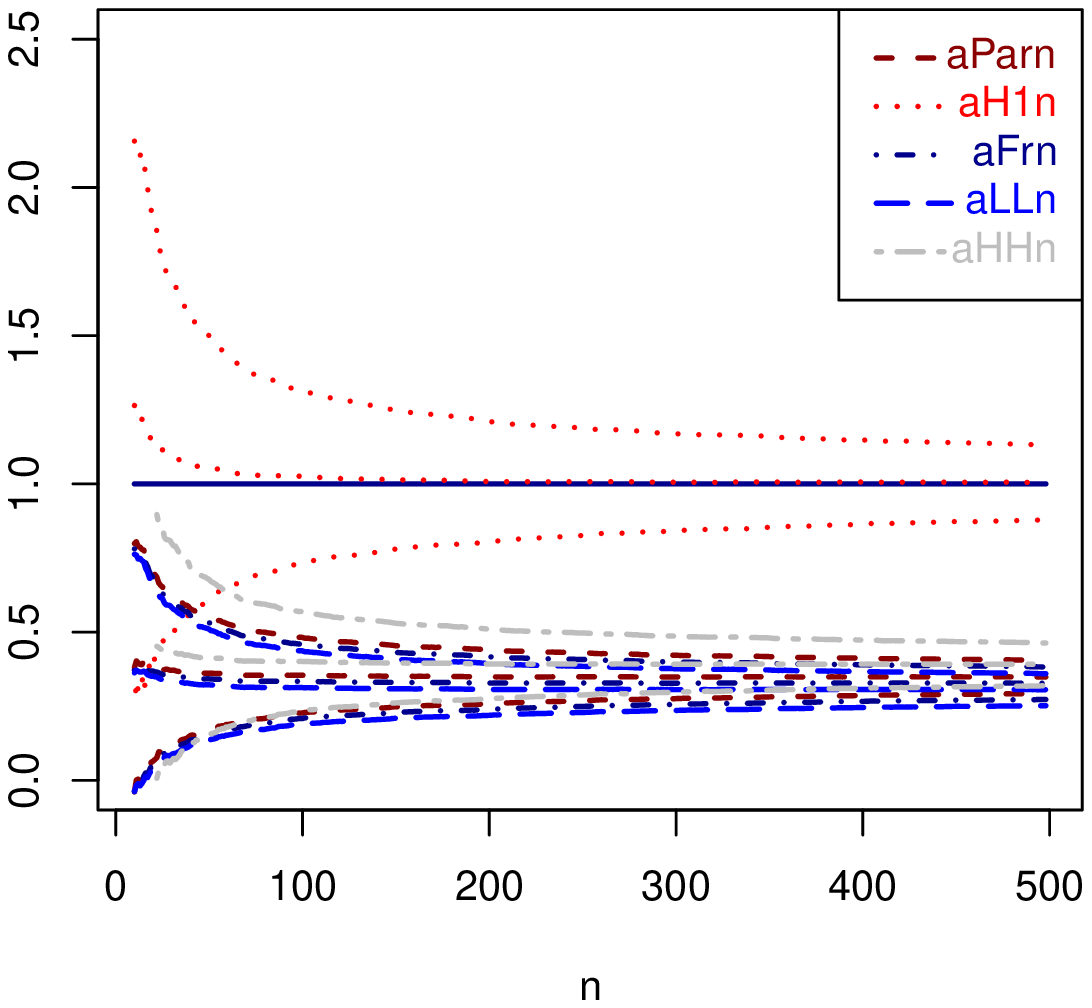}\vspace{-0.3cm}
\caption{$H_1$ case, $\alpha = 1$: Dependence of $\hat{\alpha}_{Par,n} $, $\hat{\alpha}_{Fr,n} $, $\hat{\alpha}_{HH,n}$, $\hat{\alpha}_{H_1,n}$, $\hat{\alpha}_{LL,n}$ on the sample size}
    \label{fig:SymStudy8}
\end{minipage}
\end{figure}

\begin{figure}
\begin{minipage}[t]{0.5\linewidth}
\includegraphics[scale=.47]{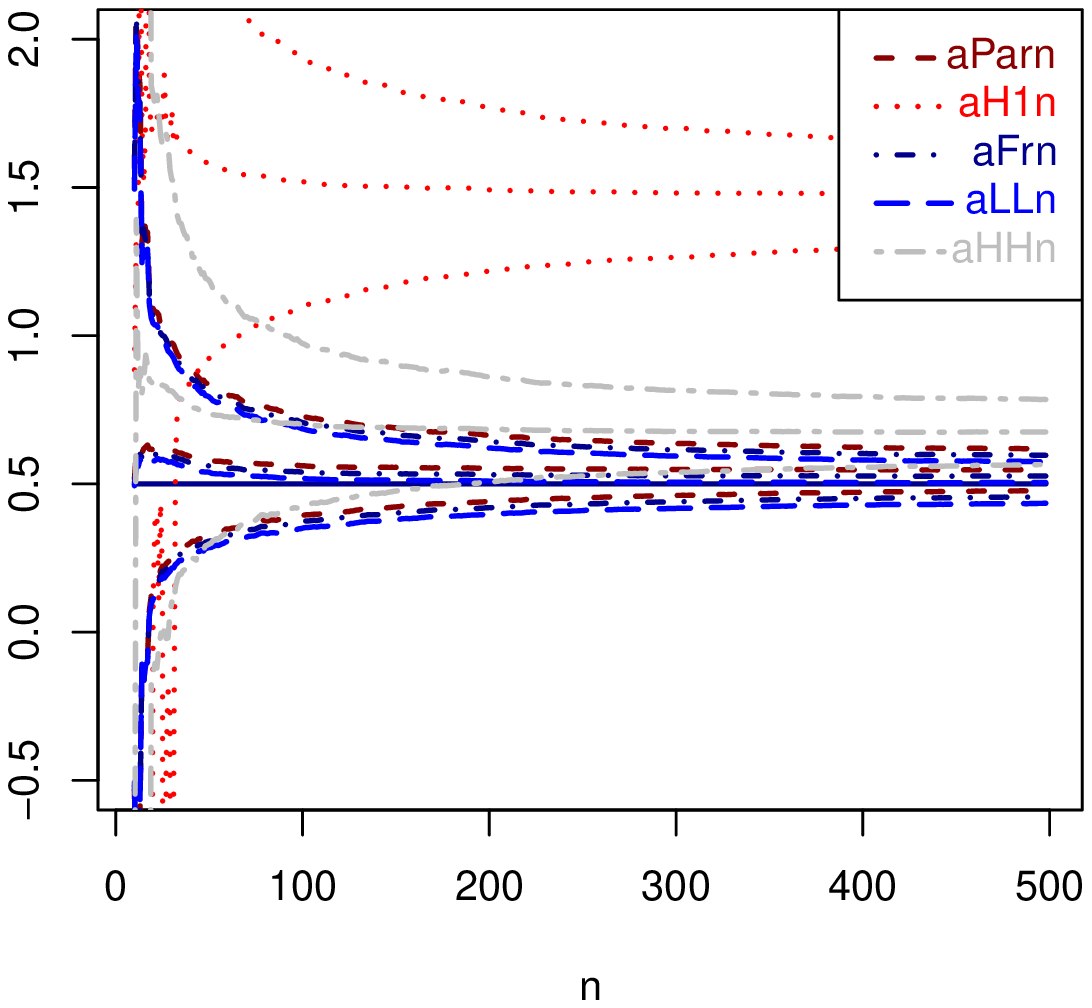}\vspace{-0.3cm}
\caption{Log-Logistic case, $\alpha = 0.5$: Dependence of $\hat{\alpha}_{Par,n} $, $\hat{\alpha}_{Fr,n} $, $\hat{\alpha}_{HH,n}$, $\hat{\alpha}_{H_1,n}$, $\hat{\alpha}_{LL,n}$ on the sample size}
    \label{fig:SymStudy9}
    \end{minipage}
\begin{minipage}[t]{0.5\linewidth}
\includegraphics[scale=.47]{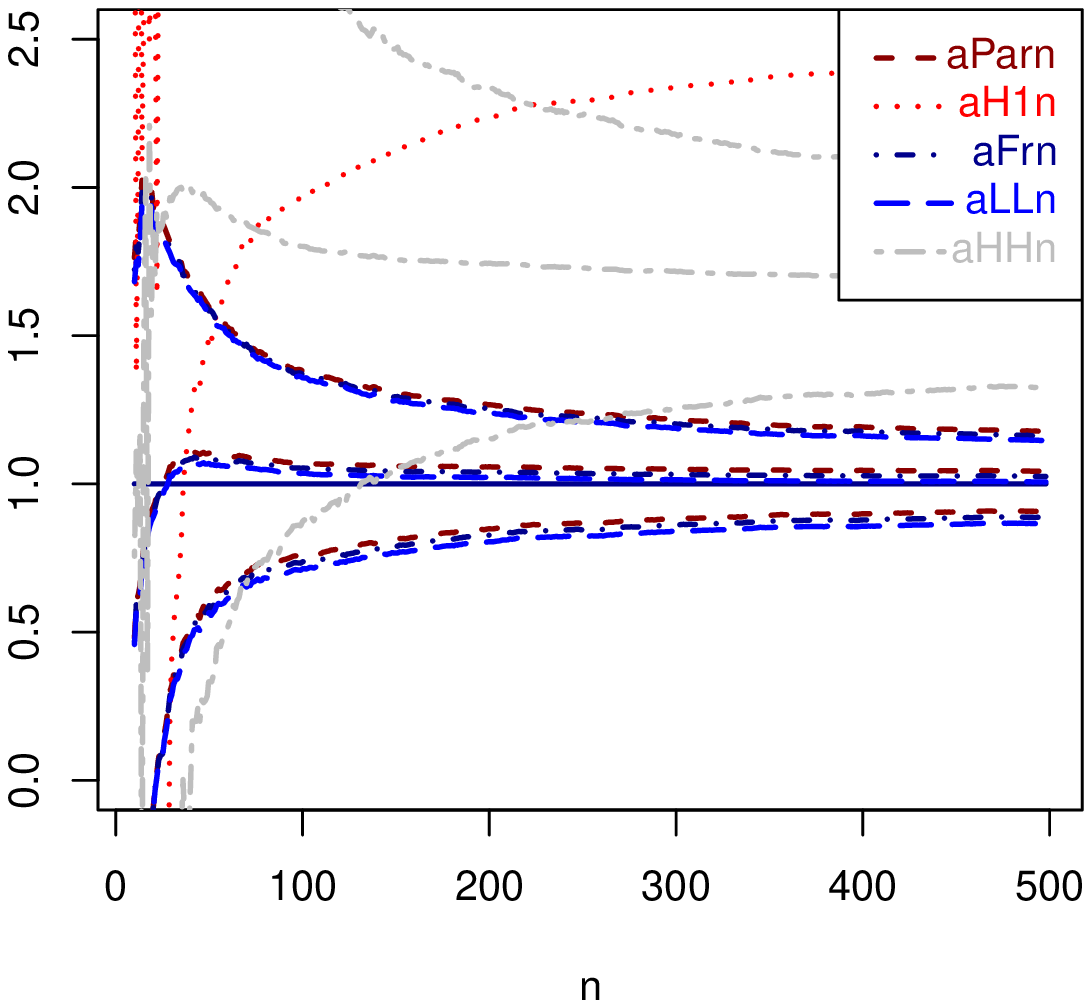}\vspace{-0.3cm}
\caption{Log-Logistic case, $\alpha = 1$: Dependence of $\hat{\alpha}_{Par,n} $, $\hat{\alpha}_{Fr,n} $, $\hat{\alpha}_{HH,n}$, $\hat{\alpha}_{H_1,n}$, $\hat{\alpha}_{LL,n}$ on the sample size}
    \label{fig:SymStudy10}
\end{minipage}
\end{figure}

The above simulation study shows that within the considered set of distributions given a small sample of observations the considered estimators outperform the properties of the well-known estimators proposed by Hill (1975) \cite{Hill}, Pickands (1975)\cite{Pickands} and Deckers-Einmahl-de Haan Dekkers (1989) \cite{Dekkers1989}. Within the right probability type, the rate of convergence of any of them increases when the sample size increases and $\alpha > 0$ decreases. However, according to our investigation, these estimators are too distribution sensitive. The biggest their advantage is that they are applicable for relatively small samples.

\section{Conclusive remarks}
\label{sec:4}

To the best knowledge of the author, a universal numerical characteristic of the tail of the c.d.f., which is invariant within distributional type (with respect to increasing affine transformation) is still not known. Here we show that probabilities of the events an observation to be $p$-outside value can be very useful in this sense. They can be used for making a reasonable classification of the tails of probability distributions. They outperform the role e.g. of the excess in characterizing the tail of the observed distribution because they do not depend on the moments of the observed r.v. and could be applied also in cases when moments do not exist. Their estimators are appropriate for usage in preliminary statistical analysis in presence of corresponding outside values. They can help the practitioners to find the most appropriate classes of probability laws for modeling the tails of the distribution of the observed r.v.  Within that family the parameter which influences the tails needs further estimation. According to our simulation study, the proposed algorithm for making estimators gives better results when $\alpha>0$ decreases. The fast rate of convergence allows one to apply these estimators also for relatively small samples. However, the main disadvantage of all these estimators is that they are distribution sensitive. The last means that their good properties may disappear if the distributional type is not correctly determined.


\section{Acknowledgements}

The author would like to thank Prof. Milan Stehlik for bringing  her to the question about statistical modeling of extremes given small samples.

\section{Apendix}

{\bf Proof of Theorem \ref{thm:monotonicity}}:
a) By definition of $L(X, p)$ and formula for derivative of the inverse function $\frac{\partial F^\leftarrow(p)}{\partial p} = \frac{1}{f(F^\leftarrow(p))}$ we obtain
$$\frac{\partial L(X, p)}{\partial p} = \frac{1}{pf[F^\leftarrow(p)]} + \frac{F^\leftarrow(1 - p) - F^\leftarrow(p)}{p^2} + \frac{1 - p }{pf[F^\leftarrow(1 - p)]}.$$
$f(x)$ is a density function of the r.v. $X$, therefore it is non-negative. The difference $F^\leftarrow(1 - p) - F^\leftarrow(p) \geq 0$ because $p \in (0; 0.5]$. Therefore $\frac{\partial L(X, p)}{\partial p} \geq 0$.

b) By definition of $R(X, p)$ and the same formula for derivative of the inverse function we have that
$$\frac{\partial R(X, p)}{\partial p} = \frac{F^\leftarrow(p) - F^\leftarrow(1 - p)}{p^2} - \frac{1}{pf[F^\leftarrow(p)]}  - \frac{1 - p }{pf[F^\leftarrow(1 - p)]}.$$
Now the difference $F^\leftarrow(1 - p) - F^\leftarrow(p) \leq 0$ because $p \in (0; 0.5]$. Therefore $\frac{\partial R(X, p)}{\partial p} \leq 0$.

c) follows by a), b), and monotonicity of probability measures.
\hfill Q.A.D.

{\bf Proof of Theorem \ref{thm:thm1}}:
b) For $c \in R$ from the definition of the quantile function we have that $F_{X+c}^\leftarrow(p) = F_X^\leftarrow(p)+c$. Therefore
\begin{eqnarray*}
 p_{R,p}(X + c) &=&  P(X+c >\frac{1}{p} F_{X+c}^\leftarrow(1-p) - \frac{1-p}{p}F_{X+c}^\leftarrow(p)) \\
&=& P(X + c > \frac{1}{p} F_{X}^\leftarrow(1-p) + \frac{c}{p} - \frac{1-p}{p}F_{X}^\leftarrow(p)- \frac{c}{p}(1-p))\\
&=& p_{R,p}(X).
\end{eqnarray*}
c) For $c > 0$ again from the definition of the quantile function $F_{cX}^\leftarrow(p) = cF_X^\leftarrow(p)$. Therefore
\begin{eqnarray*}
p_{R,p}(cX) &=& P(cX >\frac{1}{p} F_{cX}^\leftarrow(1-p) - \frac{1-p}{p}F_{cX}^\leftarrow(p)) \\
&=& P(cX > \frac{c}{p} F_{X}^\leftarrow(1-p)  - c\frac{1-p}{p}F_{X}^\leftarrow(p)) = p_{R,p}(X).
\end{eqnarray*}
d) In this case $c < 0$, therefore $F_{cX}^\leftarrow(p) = cF^\leftarrow(1-p)$ and
\begin{eqnarray*}
p_{R,p}(cX) &=& P(cX >\frac{1}{p} F_{cX}^\leftarrow(1-p) - \frac{1-p}{p}F_{cX}^\leftarrow(p)) \\
&=& P(cX > \frac{c}{p} F_{X}^\leftarrow(p)  - c\frac{1-p}{p}F_{X}^\leftarrow(1-p)) \\
&=& P(X < \frac{1}{p} F_{X}^\leftarrow(p)  - \frac{1-p}{p}F_{X}^\leftarrow(1-p)) = p_{L,p}(X).
\end{eqnarray*}
e)  Because of $g$ is a strictly increasing and continuous function we have that  $F_{g(X)}^\leftarrow(p) = g(F^\leftarrow(p))$ and $F_{g(X)}^\leftarrow(1-p) = g(F^\leftarrow(1-p))$. Therefore
\begin{eqnarray*}
  p_{R,p}(g(X)) &=& P\left[g(X) > \frac{1}{p} F_{g(X)}^\leftarrow(1-p) - \frac{1-p}{p}F_{g(X)}^\leftarrow(p)\right] \\
   &=& P\left\{g(X) > \frac{1}{p} g[F^\leftarrow(1-p)] - \frac{1-p}{p}g[F^\leftarrow(p)]\right\} \\
   &\leq& P\left\{g(X) > g[\frac{1}{p}F^\leftarrow(1-p) - \frac{1-p}{p}F^\leftarrow(p)]\right\} \\
   &=&  P\left[X > \frac{1}{p}F^\leftarrow(1-p) - \frac{1-p}{p}F^\leftarrow(p)\right] = p_{R,p}(X)
\end{eqnarray*}

These, together with the definition of $p_{R,p}(X)$ and monotonicity of probability measures entail (\ref{auxTh1e1}).

(\ref{auxTh1e2}) is a corollary of (\ref{auxTh1e3}), applied for  $\tilde{g}(x) = - g(x)$.

f) The equalities $F^\leftarrow_{g(X)}(1-p) = g(F^\leftarrow_X(p))$ and $F^\leftarrow_{g(X)}(p) = g(F^\leftarrow_X(1-p))$ entail
\begin{eqnarray*}
  p_{R, p}(g(X)) &=& P\left[g(X) > \frac{1}{p}F_{g(X)}^\leftarrow(1-p)- \frac{1-p}{p}F_{g(X)}^\leftarrow(p)\right] \\
  &=& P\left\{g(X) > \frac{1}{p} g[F^\leftarrow(p)] - \frac{1-p}{p}g[F^\leftarrow(1-p)]\right\}\\
   &\geq& P\left\{g(X) > g[\frac{1}{p} F^\leftarrow(p) - \frac{1-p}{p}F^\leftarrow(1-p)]\right\} \\
   &=& P\left[X < \frac{1}{p}F^\leftarrow(p) - \frac{1-p}{p}F^\leftarrow(1-p)\right] = p_{L,p}(X).
\end{eqnarray*}

Now the  definitions of $p_{R,p}(X)$,  and $p_{L,p}(X)$, and the monotonicity of probability measures entail (\ref{auxTh1e3}).

(\ref{auxTh1e4}) follows by (\ref{auxTh1e1}), when replace the function $g(x)$ with $\tilde{g}(x) = - g(x)$ and take into account that $p_{R, p}(-g(X)) = p_{L, p}(g(X))$.

g) As far as $F^\leftarrow_{X_{(n, n)}}(p) = F^\leftarrow_{X}(\sqrt[n]{p})$
\begin{eqnarray*}
  p_{R, p}(X_{(n, n)}) &=& P\left\{X > \frac{1}{p}F^\leftarrow_{X}(\sqrt[n]{1 - p}) - \frac{1-p}{p} F^\leftarrow_{X}(\sqrt[n]{p})\right\}\\
&=& 1 - F_X^n\left[\frac{1}{p}F^\leftarrow_{X}(\sqrt[n]{1 - p}) - \frac{1-p}{p} F^\leftarrow_{X}(\sqrt[n]{p})\right]
\end{eqnarray*}
\begin{eqnarray*}
p_{L,p}(X_{(n, n)}) &=& P\left\{X < \frac{1}{p}F^\leftarrow_{X}(\sqrt[n]{p}) - \frac{1-p}{p} F^\leftarrow_{X}(\sqrt[n]{1 - p})\right\} \\
&=&F_X^n\left[\frac{1}{p}F^\leftarrow_{X}(\sqrt[n]{p}) - \frac{1-p}{p} F^\leftarrow_{X}(\sqrt[n]{1 - p})\right]
\end{eqnarray*}
h) Using $F^\leftarrow_{X_{(1, n)}}(p) = F^\leftarrow_{X}(1-\sqrt[n]{1 - p})$ we obtain
\begin{eqnarray*}
p_{R,p}(X_{(1, n)}) &=& P\left\{X > \frac{1}{p}F^\leftarrow_{X}(1-\sqrt[n]{p}) - \frac{1-p}{p} F^\leftarrow_{X}(1-\sqrt[n]{1 - p})\right\}\\
&=&\left\{1 - F_X\left[\frac{1}{p}F^\leftarrow_{X}(1-\sqrt[n]{p}) - \frac{1-p}{p} F^\leftarrow_{X}(1-\sqrt[n]{1 - p})\right]\right\}^n
\end{eqnarray*}
\begin{eqnarray*}
p_{L,p}(X_{(1, n)}) &=& P\left\{X < \frac{1}{p}F^\leftarrow_{X}(1-\sqrt[n]{1 - p}) - \frac{1-p}{p} F^\leftarrow_{X}(1-\sqrt[n]{p})\right\} \\
&=&1 - \left\{1 - F_X \left[\frac{1}{p}F^\leftarrow_{X}(1-\sqrt[n]{1 - p}) - \frac{1-p}{p} F^\leftarrow_{X}(1-\sqrt[n]{p})\right]\right\}^n.
\end{eqnarray*}
i)  Consider $t > 0$. The relation between the quantile function of the exceedances and the c.d.f. of $F$ (see e.g. in Nair et al. (2013) \cite{Nair2013}) entails
\begin{eqnarray*}
& & p_{R,p}(X-t|X>t)\\
&=& P\left\{X - t > \frac{1}{p}F^\leftarrow_{X-t|X>t}(1 - p) - \frac{1-p}{p} F^\leftarrow_{X-t|X>t}(p)|X > t\right\} \\
&=&P\left\{X - t > \frac{1}{p}F^\leftarrow_{X}[1 - p + pF_X(p)] - \frac{t}{p} \right.\\
&-& \frac{1-p}{p} F_X^\leftarrow[p + (1-p)F_X(t)] + \left. \frac{t(1-p)}{p}|X > t\right\} \\
&=&\frac{P\left\{X > \frac{1}{p}F^\leftarrow_{X}[1 - p + pF_X(p)] - \frac{1-p}{p} F^\leftarrow_{X}[p + (1-p)F_X(t)] \right\}}{P(X>t)} \\
&=& \frac{P\left\{X > \frac{1}{p}F^\leftarrow_{X|X>t}(1 - p) - \frac{1-p}{p} F^\leftarrow_{X|X>t}(p) \right\}}{P(X>t)}\\
&=& P\left\{X > \frac{1}{p}F^\leftarrow_{X|X>t}(1 - p) - \frac{1-p}{p} F^\leftarrow_{X|X>t}(p)| X>t \right\} = p_{R,p}(X|X>t).
\end{eqnarray*}
Analogously for $p_{L,p}(X-t|X>t)$. \footnote{This property can be obtained also as a corollary of b).}

k) We obtain this property when replace the relations between the quantile functions of left-, right-, and double-truncated r.vs., $F_X^\leftarrow$ and $F_X$, i.e.
 \begin{eqnarray*}
   F_{X_{LT}}^\leftarrow(p) &=& F_X^\leftarrow(p+(1-p)F_X(l)) \\
   F_{X_{RT}}^\leftarrow(p) &=& F_X^\leftarrow(pF_X(u))   \\
   F_{X_{DT}}^\leftarrow(p) &=& F_X^\leftarrow(pF_X(u)+(1-p)F_X(l))
 \end{eqnarray*}
 in the definitions of $L(X,p)$ and $R(X,p)$. The above equalities is not difficult to calculate, and could be found e.g. in Nair et al. (2013) \cite{Nair2013}.

Finally we use the definitions of $p_{L,p}(X)$ and $p_{R, p}(X)$.

l) Assume $p \in (0, 0.5]$.

 - Case $a > 1$. In this case $F_{log_a(X)}^\leftarrow(p) = log_a[F^\leftarrow(p)]$, and the function $a^x$ is increasing in $x$, therefore
 \begin{eqnarray*}
  p_{L,p} [log_a(X)] &=& P\left[log_a(X) < \frac{1}{p}F_{log_a(X)}^\leftarrow(p) - \frac{1 - p}{p}F_{log_a(X)}^\leftarrow(1 - p)\right] \\
   &=& P\left\{log_a(X) < \frac{1}{p}log_a[F^\leftarrow(p)] - \frac{1 - p}{p}log_a[F^\leftarrow(1 - p)]\right\} \\
   &=& F \left\{\sqrt[p]{\frac{F^\leftarrow(p)}{[F^\leftarrow(1 - p)]^{1 - p}}}\right\}.
\end{eqnarray*}
\begin{eqnarray*}
  p_{R,p} [log_a(X)] &=& P\left[log_a(X) > \frac{1}{p}F_{log_a(X)}^\leftarrow(1 - p) - \frac{1 - p}{p}F_{log_a(X)}^\leftarrow(p)\right] \\
   &=& P\left\{log_a(X) > \frac{1}{p}log_a[F^\leftarrow(1 - p)] - \frac{1 - p}{p}log_a[F^\leftarrow(p)]\right\} \\
   &=& 1 - F \left\{\sqrt[p]{\frac{F^\leftarrow(1 - p)}{[F^\leftarrow(p)]^{1 - p}}}\right\}.
\end{eqnarray*}

- Case $0 < a < 1$. here we use the fact that $1/a > 1$, therefore our computations in the previous case imply
 \begin{eqnarray*}
p_{L,p} [log_a(X)] &=& p_{L,p} [-log_{1/a}(X)] = p_{R,p} [log_{1/a}(X)] \\
     &=& 1 - F \left\{\sqrt[p]{\frac{F^\leftarrow(1 - p)}{[F^\leftarrow(p)]^{1 - p}}}\right\}.
\end{eqnarray*}
\begin{eqnarray*}
  p_{R,p} [log_a(X)] &=& p_{R,p} [-log_{1/a}(X)] = p_{L,p} [log_{1/a}(X)] \\
      &=& F \left\{\sqrt[p]{\frac{F^\leftarrow(p)}{[F^\leftarrow(1 - p)]^{1 - p}}}\right\}.
\end{eqnarray*}

m)  Case $\alpha > 0$. In this case $F_{X^\alpha}^\leftarrow(p) = [F^\leftarrow(p)]^\alpha$, and the function $x^\alpha$ is increasing in $x$, therefore
 we apply Theorem 1, e) (\ref{eL}) and (\ref{eR}) and obtain the desired result.

 Case $\alpha < 0$. Now $F_{X^\alpha}^\leftarrow(p) = [F^\leftarrow(1 - p)]^\alpha$, and the function $x^\alpha$ is decreasing in $x$, therefore
 we apply Theorem 1, f) (\ref{fL}) and (\ref{fR}) complete the proof.

n) In case $\alpha \in (0, 1)$ we take into account that $F_{a^X}^\leftarrow(p) = a^[F^\leftarrow(1 - p)]$, and the functions $a^x$ and $log_a (x)$ are decreasing in $x$, then we apply Theorem 1, f) (\ref{fL}) and (\ref{fR}) and finish the proof of this case.

Analogously, if $\alpha > 1$ then $F_{a^X}^\leftarrow(p) = a^[F^\leftarrow(p)]$, and the function $x^\alpha$ is increasing in $x$. Therefore
 we apply Theorem 1, e)  and after some algebra complete the proof.

\hfill Q.A.D.

{\bf Proof of Corollary of f):} Assume $p \in (0, 0.5]$. Consider the case $P(X > 0) = 1$.  Theorem 1, f), applied for $g(x) = x^{-1}$ entails
$$p_{R, p}\left(\frac{1}{X}\right) = F\left[\frac{pF^\leftarrow(p)F^\leftarrow(1 - p)}{F^\leftarrow(1 - p) - (1 - p)F^\leftarrow(p)}\right].$$
By the definition for $p_{L, p}(X)$ and monotonicity of probability measures in order to prove that $p_{L, p}(X) \leq p_{R, p}(X^{-1})$  we need to show that
\begin{equation}\label{aux1}
\frac{1}{p}F^\leftarrow(p) - \frac{1-p}{p}F^\leftarrow(1 - p) \leq \frac{pF^\leftarrow(p)F^\leftarrow(1 - p)}{F^\leftarrow(1 - p) - (1 - p)F^\leftarrow(p)}
\end{equation}
The last inequality is equivalent to
$$1 - (1 - p)\frac{F^\leftarrow(1 - p)}{F^\leftarrow(p)} \leq p^2\frac{F^\leftarrow(1 - p)}{F^\leftarrow(p)}\left[\frac{F^\leftarrow(1 - p)}{F^\leftarrow(p)} - 1 + p\right]^{-1}.$$

For $z := \frac{F^\leftarrow(1 - p)}{F^\leftarrow(p)} \geq 1$, $z - 1 + p \geq 0$, therefore $(z - 1 + p)[1 - (1 - p)z] \leq p^2z$, which completes the proof of this part.

The assertion for the case $P(X < 0) = 1$ follows by the fact that $$p_{L, p}\left(\frac{1}{X}\right) =  p_{R, p}\left(\frac{1}{-X}\right).$$

\hfill Q.A.D.

{\bf Proof of Corollary of e):} Assume $p \in (0, 0.5]$.

 - Case $a > 1$. Because of $F_{log_a(X)}^\leftarrow(p) = log_a[F^\leftarrow(p)]$, in order to use Theorem \ref{thm:thm1}, e) for $g(x) = log_a(x)$, which is increasing in $x$ we need to prove that
$$\frac{1}{p}log_a[F^\leftarrow(1-p)] - \frac{1-p}{p}log_a[F^\leftarrow(p)] \geq log_a\left[\frac{1}{p}F^\leftarrow(1-p) - \frac{1 - p}{p}F^\leftarrow(p)\right].$$

The function $a^x$ is also increasing in $x$, therefore the above inequality is equivalent to
\begin{equation}\label{smile}
\sqrt[p]{\frac{F^\leftarrow(1 - p)}{F^\leftarrow(p)}} \geq \frac{1}{p}\frac{F^\leftarrow(1 - p)}{F^\leftarrow(p)} - \frac{1}{p} + 1.
\end{equation}
 The function $\sqrt[p]{z} - \frac{1}{p}z + \frac{1}{p} - 1$ is increasing in  $z \geq 1$, and for $p \in (0, 0.5]$, $\frac{F^\leftarrow(1 - p)}{F^\leftarrow(p)} \geq 1$. Therefore, for $z \geq 1$, $\sqrt[p]{z} - \frac{1}{p}z + \frac{1}{p} - 1 \geq 0$ proves (\ref{smile}) and completes the proof of (\ref{LogarithmsFirstinequalityR}).

By Theorem \ref{thm:thm1}, f), applied for $g(x) = log_a(x)$, (because now $log_a(x)$ is also increasing in $x$) in order to compare $p_{L, p}[log_a(X)]$ and $p_{L, p}(X)$ in (\ref{LogarithmsFirstinequalityL}) we have to show that
\begin{equation}\label{smile1}
\sqrt[p]{\frac{F^\leftarrow(p)}{[F^\leftarrow(1 - p)]^{1 - p}}} \geq \frac{1}{p}F^\leftarrow(p) - \frac{1-p}{p} F^\leftarrow(1 - p).
\end{equation}
Which is the same as
$$\sqrt[p]{\frac{F^\leftarrow(p)}{F^\leftarrow(1 - p)}} \geq \frac{1}{p}\frac{F^\leftarrow(p)}{F^\leftarrow(1 - p)} - \frac{1}{p} + 1.$$
 The function $\sqrt[p]{z} - \frac{1}{p}z + \frac{1}{p} - 1$ is decreasing in  $z \in (0, 1)$, and because of by assumption $P(X > 0) = 1$, for $p \in (0, 0.5]$, $0 < \frac{F^\leftarrow(1 - p)}{F^\leftarrow(p)} \leq 1$. Therefore, for $z \geq 1$, $\sqrt[p]{z} - \frac{1}{p}z + \frac{1}{p} - 1 \geq 0$ proves (\ref{smile1}) and (\ref{LogarithmsFirstinequalityL}).

 - Case $0 < a < 1$. Because of $1/a\geq 1$, by the previous case, and the definitions for $p_{R, p}$ and $p_{L, p}$, entail
$$p_{R, p}(log_a(X)) = p_{R, p}(-log_{1/a}(X)) = p_{L, p}(log_{1/a}(X)) \geq p_{L, p}(X)$$
  and complete the proof of (\ref{LogarithmsFirstinequalityLalessthan1}).

The fact that $1/a\geq 1$ and (\ref{LogarithmsFirstinequalityR}) entail
$$p_{L, p}(log_a(X)) = p_{L, p}(-log_{1/a}(X)) = p_{R, p}(log_{1/a}(X)) \leq p_{R, p}(X) \leq p_{R, p}(\frac{1}{a^X}), $$
and this proves inequalities in (\ref{Logarithms1}).
 \hfill Q.A.D.

{\bf Proof of Theorem \ref{thm:thm1aux}}: a) Consider $p \in (0, 0.5]$, $0 \leq \alpha_1 \leq \alpha_2$ and $P$ almost sure positive r.v. $X$, i.e.  $P(X > 0) = 1$. According to monotonicity of probability measures, the definition of $p_{R,p}$, and the equalities $F_{X^{\alpha_i}}^\leftarrow(p) = [F^\leftarrow(p)]^{\alpha_i}$, $i = 1, 2$  we need to show that
\begin{equation}\label{aux2}
\sqrt[\alpha_1]{\frac{1}{p}[F^\leftarrow(1-p)]^{\alpha_1} - \frac{1 - p}{p}[F^\leftarrow(p)]^{\alpha_1}} \geq \sqrt[\alpha_2]{\frac{1}{p}[F^\leftarrow(1-p)]^{\alpha_2} - \frac{1 - p}{p}[F^\leftarrow(p)]^{\alpha_2}}.
\end{equation}
It is the same as
$$\sqrt[\alpha_1]{\frac{1}{p}\left[\frac{F^\leftarrow(1-p)}{F^\leftarrow(p)}\right]^{\alpha_1} - \frac{1}{p} + 1} \geq \sqrt[\alpha_2]{\frac{1}{p}\left[\frac{F^\leftarrow(1-p)}{F^\leftarrow(p)}\right]^{\alpha_2} - \frac{1}{p} + 1}.$$

Denote by $z:=\frac{F^\leftarrow(1-p)}{F^\leftarrow(p)} \geq 1$. The last inequality is true, because of for any fixed $z \geq 1$ and $p \in (0; 0.5]$ the function
$$h(\alpha) = \sqrt[\alpha]{\frac{1}{p}z^{\alpha} - \frac{1}{p} + 1}$$
is decreasing in $\alpha > 0$. Therefore (\ref{aux2}) is also true, and the proof of (\ref{Cor2ofe1}) is completed.

b) The r.v. $X$ is almost sure positive, so we have the same for $X^\alpha$. Therefore, by the definition of $p_{L, p}(X^\alpha)$ it is enough to show that, in this case $L(X^\alpha, p) \leq 0$. Now we use the equality $F^\leftarrow_{X^\alpha}(p) = [F^\leftarrow_{X}(p)]^\alpha$ in the definition of $L(X^\alpha, p)$ and obtain that given the condition in c), the value of
$$L(X^\alpha, p) = \frac{1}{p}\left\{[F^\leftarrow(p)]^\alpha - (1 - p)[F^\leftarrow(1 - p)]^\alpha\right\} \leq 0.$$

c) It is enough to prove the second inequality in (\ref{Cor2ofe2}). It is equivalent to
$$P\left[X < \frac{1}{p}F^\leftarrow(p) - \frac{1 - p}{p}F^\leftarrow(1 - p)\right] \geq P\left[X^\alpha < \frac{1}{p}F_{X^\alpha}^\leftarrow(p) - \frac{1 - p}{p}F_{X^\alpha}^\leftarrow(1 - p)\right],$$
and using the equality $F_{X^\alpha}^\leftarrow(p) = [F^\leftarrow(p)]^\alpha$ it is the same as
$$P\left[X < \frac{1}{p}F^\leftarrow(p) - \frac{1 - p}{p}F^\leftarrow(1 - p)\right] \geq P\left[X <\sqrt[\alpha]{\frac{1}{p}[F^\leftarrow(p)]^\alpha - \frac{1 - p}{p}[F^\leftarrow(1 - p)]^\alpha}\right].$$

The monotonicity of probability measures entails that the above inequality would be true if
$$\frac{1}{p}F^\leftarrow(p) - \frac{1 - p}{p}F^\leftarrow(1 - p) \geq \sqrt[\alpha]{\frac{1}{p}[F^\leftarrow(p)]^\alpha - \frac{1 - p}{p}[F^\leftarrow(1 - p)]^\alpha},\,\, i.e.$$
$$\left[\frac{1}{p}F^\leftarrow(p) - \frac{1 - p}{p}F^\leftarrow(1 - p)\right]^\alpha \geq \frac{1}{p}[F^\leftarrow(p)]^\alpha - \frac{1 - p}{p}[F^\leftarrow(1 - p)]^\alpha, \,\, that \,\, is$$
\begin{equation}\label{aux3}
\left[\frac{1}{p} - \frac{1 - p}{p}\frac{F^\leftarrow(1 - p)}{F^\leftarrow(p)}\right]^\alpha \geq \frac{1}{p}- \frac{1 - p}{p}\left[\frac{F^\leftarrow(1 - p)}{F^\leftarrow(p)}\right]^\alpha.
\end{equation}
Again denote by $z:=\frac{F^\leftarrow(1-p)}{F^\leftarrow(p)} \geq 1$, and consider the function
$$t(z): = \frac{1}{p} - \frac{1 - p}{p}z^\alpha - \left[\frac{1}{p} - \frac{1 - p}{p}z\right]^\alpha.$$
Given $\alpha > 1$, it is decreasing for $z \geq 1$ and $t(1) = 0$. Therefore $t(z) \leq 0$, for $z \geq 1$.  This entails (\ref{aux3}) and completes the proof of (\ref{Cor2ofe2}).   \hfill Q.A.D.

{\bf Proof of Theorem \ref{thm:ThmGeneralAsumptoticNormalityFences}:} Let us fix $s \in N$. From the theorem about the joint asymptotic normality of the order statistics, because the limit exists we have that for any subsequence $\{n: n = (s+1)k-1, s \in N\}$,
 $$\sqrt{(s+1)k-1}\left(\begin{array}{c}
                             X_{(k, (s+1)k-1)} - F^\leftarrow(\frac{1}{s+1}) \\
                             X_{(ki, (s+1)k-1)} - F^\leftarrow(\frac{i}{s+1}) \\
                           \end{array}
                         \right) \stackrel{d}{\to} N\left[\left(\begin{array}{c}
                         0\\
                         0 \\
                           \end{array}\right); D\right], \quad k \to \infty,$$
 where the asymptotic covariance matrix of this bivariate distribution is
 $$D = \frac{1}{(s+1)^2}\left(\begin{array}{cc}
 \frac{s}{f^2\left[F^\leftarrow(\frac{1}{s+1})\right]} & \frac{1}{f\left[F^\leftarrow(\frac{1}{s+1})\right]f\left[F^\leftarrow(\frac{s}{s+1})\right]} \\
 \frac{1}{f\left[F^\leftarrow(\frac{1}{s+1})\right]f\left[F^\leftarrow(\frac{s}{s+1})\right]} & \frac{s}{f^2\left[F^\leftarrow(\frac{s}{s+1})\right]} \\
 \end{array}\right)$$
 and the asymptotic correlation between these two order statistics is $\frac{1}{i}$.

 1.) Consider the function $g(x, y) = (s+1)x - sy$. For $x > 0$ and $y > 0$ it is continuously differentiable.
 The asymptotic mean is
 \begin{eqnarray*}
  \lim_{k \to \infty} E L_n\left(X, \frac{1}{s+1}\right) &=& g\left[F^\leftarrow(\frac{1}{s+1}), F^\leftarrow(\frac{s}{s+1})\right]\\
    &=& (s+1)F^\leftarrow(\frac{1}{s+1}) - i F^\leftarrow(\frac{s}{s+1}) = L\left(X,\frac{1}{s+1}\right).
 \end{eqnarray*}

The Jacobian of the transformation is
$$I_L : = \left[\frac{\partial g(x, y)}{\partial x}, \frac{\partial g(x, y)}{\partial y}\right] = \left(s+1, -s\right).$$
Now we apply the Multivariate Delta method Sobel (1982) \cite{MultivariateDeltaMethod}, and obtain that the asymptotic variance of $L_n\left(X, \frac{1}{s+1}\right)$ is
 \begin{eqnarray*}
  V_L : &=& I_L \times D \times I_L' = \left(s+1, -s\right)\left|_{x = F^\leftarrow(\frac{1}{s+1}), y = F^\leftarrow(\frac{s}{s+1})}\right. \\
    &\times&  \frac{1}{(s+1)^2}\left(\begin{array}{cc}
 \frac{s}{f^2\left[F^\leftarrow(\frac{1}{s+1})\right]} & \frac{1}{f\left[F^\leftarrow(\frac{1}{s+1})\right]f\left[F^\leftarrow(\frac{s}{s+1})\right]} \\
 \frac{1}{f\left[F^\leftarrow(\frac{1}{s+1})\right]f\left[F^\leftarrow(\frac{s}{s+1})\right]} & \frac{i}{f^2\left[F^\leftarrow(\frac{s}{s+1})\right]} \\
 \end{array}\right) \\
    &\times&  \left(
                             \begin{array}{c}
                               s+1 \\
                               -s \\
                             \end{array}
                           \right)\left|_{x = F^\leftarrow(\frac{1}{s+1}), y = F^\leftarrow(\frac{s}{s+1})}\right. \\
   &=& \frac{1}{(s+1)^2} \\
 &\times& \left[\left(s+1, -s\right)\left(\begin{array}{cc}
 \frac{i}{f^2(x)} & \frac{1}{f(x)f(y)} \\
 \frac{1}{f(x)f(y)} & \frac{i}{f^2(y)} \\
 \end{array}\right) \left(
                             \begin{array}{c}
                               s+1 \\
                               -s\\
                             \end{array}
                           \right)\right]\left|_{x = F^\leftarrow(\frac{1}{s+1}), y = F^\leftarrow(\frac{s}{s+1})}\right. \\
    &=& \frac{1}{(s+1)^2} \left[\frac{s(s+1)^2}{f^2(x)} - \frac{2s(s+1)}{f(x)f(y)} + \frac{s^3}{f^2(y)}\right] \left|_{x = F^\leftarrow(\frac{1}{s+1}), y = F^\leftarrow(\frac{s}{s+1})}\right. \\
    &=& \frac{s}{(s+1)^2} \left[\frac{(s+1)^2}{c_{F,s}^2} - \frac{2(s+1)}{c_{F,s}d_{F,s}} + \frac{s^2}{d_{F,s}^2}\right].
 \end{eqnarray*}

2.) Analogously, because of
\begin{eqnarray*}
 R_{n}\left(F, \frac{1}{s + 1}\right) &=& (s + 1)X_{(ks, (s + 1)k-1)} - i X_{(k, (s + 1)k-1)}, \,\, and  \\
 R\left(X, \frac{1}{s+1}\right)  &=& (s + 1)F^\leftarrow\left(\frac{s}{s + 1}\right) + s F^\leftarrow\left(\frac{1}{s + 1}\right),
\end{eqnarray*}
we consider function $g_R(x, y) = -sx + (s + 1)y$.

For $x > 0$ and $y > 0$ it is continuously differentiable.

The asymptotic mean is
  \begin{eqnarray*}
  \lim_{k \to \infty} E R_{n}\left(F, \frac{1}{s + 1}\right) &=& g_R\left[F^\leftarrow(\frac{1}{s+1}), F^\leftarrow(\frac{s}{s+1})\right]\\
   &=& -sF^\leftarrow(\frac{1}{s+1}) + (s + 1)F^\leftarrow(\frac{s}{s+1}) = R\left(X, \frac{1}{s+1}\right).
 \end{eqnarray*}

In order to obtain the asymptotic variance of $R_n\left(X, \frac{1}{s+1}\right)$ we calculate the Jacobian of the transformation. It is
$$I_R : = \left[\frac{\partial g_R(x, y)}{\partial x}, \frac{\partial g_R(x, y)}{\partial y}\right] = \left(s+1, -s\right).$$

Now we apply the Multivariate Delta method (see e.g. Sobel (1982) \cite{MultivariateDeltaMethod}), calculate $I_R \times D \times I_R'$  and obtain the asymptotic variance of  $R_{n}\left(F, \frac{1}{s + 1}\right)$ is
$$V_R := I_R \times D \times I_R' = \frac{s}{(s+1)^2} \left[\frac{s^2}{c_{F,s}^2} - \frac{2(s+1)}{c_{F,s}d_{F,s}} + \frac{(s + 1)^2}{d_{F,s}^2}\right].$$
 \hfill Q.A.D.









\begin{thebibliography}{}

\bibitem{Cordeiro} Alexander, C., Cordeiro, G. M., Ortega, E. M.,  Sarabia, J. M.: Generalized beta-generated distributions. Computational Statistics and Data Analysis. 56(6), (2012) 1880--1897.

\bibitem{Arnold1992}  Arnold, B. C., Balakrishnan, N.,  Nagaraja,  H. N.: A first course in order statistics. In Applied Mathematics SIAM, 54,  Philadelphia (1992).

\bibitem{Arnold2015}  Arnold, B. C.: Pareto distributions, Second Edition, Chapman and Hall$/$CRC Press Taylor \& Francis Group, Boca Raton, London, New York (2015).

\bibitem{BalandaMacGillary1988} Balanda, K. P., MacGillivray, H. L.: Kurtosis: a critical review. The American Statistician. 42(2), (1988) 111--119.

\bibitem{BalandaMacGillary} Balanda, K. P., MacGillivray, H. L.:  Kurtosis and spread.  Canadian Journal of Statistics. 18(1), (1990) 17--30.

\bibitem{Beran} Beran, J., Dieter Schell, Stehlik, M.: The harmonic moment tail index estimator: asymptotic distribution and robustness. Annals of the Institute of Statistical Mathematics.  66(1), (2014) 193--220.

\bibitem{Burr} Burr, I. W.: Cumulative frequency functions. Annals of Mathematical Statistics. 13 (2), (1942) 215-–232.

\bibitem{CaeiroGomes} Caeiro, F., Gomes, M. I., Beirlant, J.,  de Wet, T.: Mean-of-order p reduced-bias extreme value index estimation under a third-order framework. Extremes. 19(4), (2016) 561--589.

\bibitem{Cadwell} Cadwell, J.H.: The distribution of quasi-ranges in samples from a normal population. The Annals of Mathematical Statistics. 24(4), (1953) 603--613.

\bibitem{Chu} Chu, J., T.: Some uses of quasi-ranges. The Annals of Mathematical Statistics. 28(1), (1957) 173--180.

\bibitem{CormannAndReiss2009} Cormann, U., Reiss, R.-D.: Generalizing the Pareto model to the log Pareto model and statistical inference. Extremes. 12, (2009) 93–-105.

\bibitem{Daouia2018a} Daouia, A., Girard, S.,  Stupfler, G.: Estimation of tail risk based on extreme expectiles. Journal of the Royal Statistical Society: Series B (Statistical Methodology). 80(2), (2018) 263--292.

\bibitem{deHaanFerreira} De Haan, L., Ferreira, A.: Extreme Value Theory: An introduction. Springer Series in Operations Research and Financial Engineering. Springer Science and Business Media, New York (2006).

\bibitem{deHaanStadtmueller} De Haan, L., Stadtm\"uller, U.: Generalized regular variation of second order. Journal of the Australian Mathematical Society, 61(3),  (1996) 381--395.

\bibitem{Dekkers1989}  Dekkers, A.L.,  Einmahl, J.H.,  De Haan, L.: A moment estimator for the index of an extreme-value distribution. The Annals of Statistics, (1989) 1833-–1855.

\bibitem{Dembinska2017Metrika} Dembi$\acute{n}$ska, A.: An ergodic theorem for proportions of observations that fall into random sets determined by sample quantiles. Metrika. 80(3), (2017) 319--332.

\bibitem{Dembinska2012NZ} Dembi$\acute{n}$ska, A.: Limit theorems for proportions of observations falling into random regions determined by order statistics. Australian and New Zealand Journal of Statistics. 54(2), (2012) 199--210.

\bibitem{Devore} Devore, J.L.: Probability and Statistics for Engineering and the Sciences. Cengage Learning, Australia, Brazil, Mexico, Singapur, United Kindom, United States (2015).

\bibitem{EinmahlGuillou} Einmahl, J.H.J., Fils-Villetard, A., Guillou, A.: Statistics of extremes under random censoring. Bernoulli. 14(1), (2008) 207--227.

\bibitem{EMK}  Embrechts, P., Kl\"uppelberg, C., Mikosch, T.: Modelling Extremal Events
for Insurance and Finance. Springer-Verlag, Berlin, Heidelberg, New York, London, Paris, Tokio, Hong Kong, Barcelona, Budapest (1997).

\bibitem{Eugene} Eugene, N., Lee, C., Famoye, F.:  Beta-normal distribution and its applications. Communications in Statistics – Theory and Methods. 31, (2002)  497--512.

\bibitem{Falk2004} Falk, M., Huesler, J., Reiss, R-D.: Functional Laws of Small Numbers. In: Falk, M.,  Huesler, J., Reiss, R-D. (eds.)  Laws of Small Numbers. Extremes and Rare Events, pp. 3--22, Springer Science and Business Media, Basel Birkh$\ddot{a}$user Verlag AG, Basel AG (2004).

\bibitem{FisherTippet} Fisher, R.A., Tippett, L.H.C.: Limiting forms of the frequency distribution of the largest or smallest member of a sample. Mathematical Proceedings of the Cambridge Philosophical Society. Cambridge University Press. 24(2), (1928) 180--190.

\bibitem{GelukdeHaan1997} Geluk, J., de Haan, L.,   Resnick,S.I., Staarica, C.:   Second-order regular variation, convolution and the central limit theorem,  Stochastic Process. Appl. 69(2), (1997)  139--159.

\bibitem{Gnedenko1943} Gnedenko, B.V.: On increments of homogeneous random processes with independent increments. Proceedings of the Russian Academy of Sciences. Series Mathematics, 7(2), (1943)  89-110 (in russian).

\bibitem{Gumbel1944} Gumbel, E.J.: Ranges and midranges. The Annals of Mathematical Statistics. 15(4), (1944) 414--422.

\bibitem{Gumbel1958} Gumbel, E.J.: Statistics of extremes. Columbia University press, New York 247 (1958).

\bibitem{GuoFeng2013} Guo, Bai-Ni, Feng Qi: Refinements of lower bounds for polygamma functions. Proceedings of the American Mathematical Society, (2013) 1007--1015.

\bibitem{Harter} Harter, H.L.: The use of sample Quasi-ranges in estimating population standard deviation. The Annals of Mathematical Statistics. 30(4), (1959) 980--999.

\bibitem{Hill} Hill, B.:  A simple general approach to inference about the tail of a distribution, The Annals of Statistics, 3(5), 1163--1174 (1975)

\bibitem{Huisman} Huisman, R.,  Koedi, K.G., Kool, C.J.M., Franz, P.:  Tail-index estimates in small samples, Journal of Business \& Economic Statistics, 19(2), (2001) 208--216.

\bibitem{Hyndman} Hyndman, R.J., Yanan, F.: Sample quantiles in statistical packages. The American Statistician.  50(4), (1996), 361--365.

\bibitem{OurExtremes} Jordanova, P.K., Fabian, Z., Hermann, P., Strelec, L., Rivera, A., Girard, S., Tores S., Stehlik, M.: Weak properties and robustness of t-Hill estimators. Extremes. 19(4), (2016) 591-–626.

\bibitem{JordanovaMilan2012} Jordanova P.K., Stehlik, M., Fabian Zd., Strelec L.: On Estimation And Testing For Pareto Tails, Pliska Studia Mathematica Bulgarica. 22, (2013) 89--108.

\bibitem{MoniPoli2017} Jordanova, P.K., Petkova, M.P.:  Measuring heavy-tailedness of distributions, AIP Conference Proceedings. AIP Publishing. 1910(1), (2017) 060002-1--060002-8.

\bibitem{MoniPoli2018} Jordanova, P.K., Petkova, M.P.:  Tails and probabilities for extreme outliers. Application of Mathematics in Technical and Natural Sciences. AIP Conference Prococeedings. USA: AIP Publishing LLC. 2025, (2018) 030002-1–-030002-9.

\bibitem{Langford}  Langford, E.:  Quartiles in elementary statistics. Journal of Statistics Education.  14(3), (2006). Available online on http://jse.amstat.org/v14n3/langford.html,  Accessed 10 January 2019.

\bibitem{MarinelliCarlo2007} Marinelli, C., d\'Addona, S., Rachev, S., T.: A comparison of some univariate models for value-at-risk and expected shortfall. International Journal of Theoretical and Applied Finance. 10(06), (2007) 1043--1075.

\bibitem{Monsteller1946} Monsteller, F.: On some useful \'inefficient\' statistiscs, Annals of Mathematical Statistics. 17, (1946) 377--408.

\bibitem{Nair2013} Nair, N.U., Sankaran, P.G.,  Balakrishnan, N.:  Quantile-based reliability analysis. Birkh$\ddot{a}$user, Basel (2013).

\bibitem{ClNeves}  Neves, Cl., Fraga Alves, M.I.: Ratio of Maximum to the Sum for Testing Super Heavy Tails. In: Minguez, R.,  Sarabia, J.-M.,  Balakrishnan, N., Arnold, B.C. (eds.) Advances in Mathematical and Statistical Modeling. 181--194. Birkhauser, Boston (2008).

\bibitem{Nevzorov} Nevzorov, V.B.: Records: mathematical theory. American Mathematical Society, United States (2001).

\bibitem{Pancheva}  Pancheva, E.: Limit theorems for extreme order statistics under nonlinear normalization. Lecture Notes in Mathematics. 1155, (1985) 284--309.

\bibitem{JordanovaPancheva}  Pancheva, E., Jordanova, P.:  Weak Asymptotic Results for t-Hill Estimator.  Comptes Rend. Acad. Bulg. Sci.  65(12), (2012) 1649-1656.

\bibitem{Parzen}  Parzen, E.: Nonparametric statistical data modeling. Journal of the American statistical association. 74(365), (1979) 105--121.

\bibitem{Paulauskas} Paulauskas, V., Vaiciulis, M.: A class of new tail index estimators. Annals of the Institute of Statistical Mathematics, 69(2), (2017) 461-487.

\bibitem{Pickands}  Pickands, J. III:  Statistical inference using extreme order statistics.  The Annals of Statistics. 3(1), (1975)  119–-131.

\bibitem{R} R Development Core Team: R: A Language and Environment for Statistical Computing, R Foundation for Statistical Computing, (2018),  http://www.R-project.org, Accessed 10 January 2019.

\bibitem{Ravi}  Ravi, S., Saeb, A.:  A note on entropies of l-max stable, p-max stable, generalized Pareto and generalized log-Pareto distributions. ProbStat Forum. 5 July, (2012) 62–-79.

\bibitem{Rider} Rider, P.R.: Quasi-ranges of samples from an exponential population. The Annals of Mathematical Statistics. 30(1), (1959) 252--254.

\bibitem{Resnick87} Resnick, S.I.:  Extreme Values, Regular Variation and Point Processes, Springer Series in Operations Research and Financial Engineering,  Springer-Verlag, New York (1987).

\bibitem{SarhanGeneral} Sarhan, A.E.: Estimation of the mean and standard deviation by order statistics. The Annals of Mathematical Statistics. 25(2),  (1954) 317--328.

\bibitem{SarhanExponential} Sarhan, A.E.:  Greenberg, B. G., Ogawa, J.: Simplified estimates for the exponential distribution. The Annals of Mathematical Statistics. 34(1), (1963) 102-116.

\bibitem{Smirnov1949} Smirnov, N.V.: Limit distributions for the terms of a variational series. Trudy Matematicheskogo Instituta imeni VA Steklova. 25, (1949) 3--60.

\bibitem{MultivariateDeltaMethod} Sobel, M.E.: Asymptotic confidence intervals for indirect effects in structural equation models. Sociological methodology. 13, (1982) 290-312.

\bibitem{JordanovaStehlik2018} Soza, L.N., Jordanova, P., Nicolis, O., Strelec, L., Stehlik, M.:  Small sample robust approach to outliers and correlation of Atmospheric Pollution and Health Effects in Santiago de Chile. Chemometrics and Intelligent Laboratory Systems. 185, (2019) 73--84.

\bibitem{Stehlik2010} Stehlik M.,  Potocky, R.,  Waldl, H.,  Fabian, Z.:  On the favorable estimation
for fitting heavy tailed data. Computational Statistics. 25 (3), (2010) 485–-503.

\bibitem{Tukey1977}  Tukey, J.W.: Exploratory data analysis, Addison-Wesley, Reading  (1977).

\bibitem{Wilks1948} Wilks, S. S. Order statistics. Bull. Amer. Math. Soc. 54, no. Number 1, Part 1, (1948) 6--50. https://projecteuclid.org/euclid.bams/1183511502

\end{thebibliography}

\end{document}